\documentclass{amsart}
\usepackage{amsfonts}

\newtheorem{thm}{Theorem}

\newtheorem{lem}[thm]{Lemma}
\newtheorem{cor}[thm]{Corollary}

\newtheorem{prop}[thm]{Proposition}

\newtheorem{conj}[thm]{Conjecture}
\newtheorem{defn-thm}[thm]{Definition--Theorem}  
\newtheorem{defn-prop}[thm]{Definition--Proposition}  

\theoremstyle{definition}
\newtheorem{defn}[thm]{Definition}

\newtheorem{exmp}[thm]{Example}

\newtheorem{ques}[thm]{Question}    

\newtheorem{rem}[thm]{Remark}
\newtheorem{ack}{Acknowledgments\!\!\!}

\theoremstyle{remark}

\renewcommand{\c}[0]{{\mathbb C}}

\newcommand{\be}{\begin{equation}}
\newcommand{\ea}{\end{array}}\newcommand{\ee}[1]{\label{#1}\end{equation}}

\newcommand{\qtq}[1]{\quad\mbox{#1}\quad}

\def\fract#1#2{\raise4pt\hbox{$ #1 \atop #2 $}}
\def\decdnar#1{\phantom{\hbox{$\scriptstyle{#1}$}}
\left\downarrow\vbox{\vskip15pt\hbox{$\scriptstyle{#1}$}}\right.}

\def\bfa{{\bf a}}

\def\bfw{{\bf w}}
\def\bfx{{\bf x}}

\def\bfz{{\bf z}}

\def\calc{{\mathcal C}}
\def\calo{{\mathcal O}}
\def\calu{{\mathcal U}}
\def\cald{{\mathcal D}}
\def\cale{{\mathcal E}}
\def\calf{{\mathcal F}}

\def\calj{{\mathcal J}}
\def\calk{{\mathcal K}}
\def\call{{\mathcal L}}

\def\calo{{\mathcal O}}

\def\cals{{\mathcal S}}

\def\calu{{\mathcal U}}

\def\calz{{\mathcal Z}}

\def\bbc{{\mathbb C}}

\def\bbh{{\mathbb H}}
\def\bbi{{\mathbb I}}

\def\bbp{{\mathbb P}}

\def\bbr{{\mathbb R}}

\def\bbz{{\mathbb Z}}

\def\gra{\alpha}

\def\grd{\delta}

\def\grg{\gamma}

\def\grl{\lambda}

\def\gro{\omega}

\def\grr{\rho}

\def\grt{\tau}

\def\grD{\Delta}

\def\grG{\Gamma}

\def\grL{\Lambda}
\def\grO{\Omega}

\def\gsp1{{\mathfrak s}{\mathfrak p}(1)}

\def\ker{\hbox{ker}}

\def\barj{\bar{j}}

\def\gh{{\mathfrak h}}

\def\la#1{\hbox to #1pc{\leftarrowfill}}
\def\ra#1{\hbox to #1pc{\rightarrowfill}}

\def\Ke{K\"ahler-Einstein }

\begin{document}
\bibliographystyle{amsalpha}

\title{Sasakian Geometry, Hypersurface Singularities, and Einstein Metrics}
\author{Charles P. Boyer and Krzysztof Galicki}
\address{Department of Mathematics and Statistics,
University of New Mexico,
Albuquerque, NM 87131.}
\address{Max-Planck-Institut f\"ur Mathematik, D53111 Bonn, Germany; on leave from
Department of Mathematics and Statistics,
University of New Mexico,
Albuquerque, NM 87131}

\email{cboyer@math.unm.edu}
\email{galicki@math.unm.edu}

\maketitle
\bigskip
\section{Introduction}
\bigskip
This review article has grown out of notes for the three lectures
the second author presented during the XXIV-th Winter School of {\it Geometry and
Physics} in Srni, Czech Republic, in January of 2004. Our purpose is twofold. We want give a 
brief
introduction to some of the techniques we have developed
over the last 5 years while, at the same time, we
summarize all the known results. We do not give any technical details other than
what is necessary for the clarity of the exposition. In conclusion we
would like to argue that
Sasakian geometry has emerged as one of the most powerful
tools of constructing and proving existence of special Riemannian metrics,
such as Einstein metrics or metrics of positive Ricci curvature, on a wide range
of odd-dimensional manifolds. The key geometric object in the theory is
that of a contact structure (hence, only odd dimensions) together with a
Riemannian metric naturally adapted to the contact form. Sasakian metrics in
contact geometry are analogous to the K\"ahler metrics in the symplectic
case.

We begin with basic facts about contact and Sasakian manifolds
after which we focus on exploring  the fundamental relation
between the Sasakian and the transverse K\"ahler geometry. In this
context positive Sasakian, Sasakian-Einstein, and 3-Sasakian
manifolds are introduced. In Section 3 we present all known
constructions of 3-Sasakian manifolds. These come as V-bundles
over compact quaternion K\"ahler orbifolds and large families can
be explicitly obtained using symmetry reduction. We also discuss
Sasakian-Einstein manifolds which are not 3-Sasakian. Here there
has been only one effective method of producing examples, namely
by representing a Sasakian-Einstein manifold as the total space of
an $S^1$ Seifert bundle over a K\"ahler-Einstein orbifold. In the
smooth case with a trivial orbifold structure, this construction
goes back to Kobayashi \cite{Kob56}. Any smooth Fano variety
$\calz$ which admits a K\"ahler-Einstein metric can be used for
the base of a unique simply connected circle bundle $P$ which is
Sasakian-Einstein. Any Sasakian-Einstein manifold obtained this
way is automatically regular. It is clear that, in order to get
non-regular examples of Sasakian-Einstein structures, one should
replace the smooth Fano structure with a Fano orbifold. This was
done in \cite{BG00a} where we generalized the Kobayashi
construction to V-bundles over Fano orbifolds. However, at that
time, with the exception of twistor spaces of known 3-Sasakian
metrics,  compact Fano orbifolds known to admit orbifold
K\"ahler-Einstein metrics were rare. The first examples of non-regular
Sasakian-Einstein manifolds which are not 3-Sasakian were obtained
in \cite{BG00a}. There we observed that Sasakian-Einstein
manifolds have the structure of a monoid under a certain ``join"
operation. A join $M_1\star M_2$ of non-regular 3-Sasakian
manifold and a regular Sasakian-Einstein manifold (say an
odd-dimensional sphere) is automatically a non-regular
Sasakian-Einstein space. The problem is that our join construction
produced new examples starting in dimension 9 and higher. It gave
nothing new in dimensions 5 and 7. The construction of
5-dimensional examples followed, however, a year later, and with
that, new non-regular Sasakian-Einstein examples in every dimension
could be obtained.

In  \cite{BG01b} the authors constructed three 5-dimensional examples using
the results of Demailly and Koll\'ar \cite{DeKo01}.
At the same time it became clear that the so-called continuity method
as being applied by Demailly and Koll\'ar to Fano orbifolds gave a cornucopia
of new examples of K\"ahler-Einstein orbifolds \cite{Ara02, BG03, BG03p, BGN02a, BGN02b,
BGN03c, BGK03, BGKT03, jk1, jk2, Kol04}.
A particular illustration
of how the method works comes from one example in classical differential topology.
It is well-known \cite{Tak,BG01b} that any link of isolated hypersurface singularities has a
natural Sasakian structure. The transverse K\"ahler geometry in such a situation
is induced by a K\"ahlerian embedding of  a complex hypersurface in
an appropriate weighted projective space. Section 6 describes Sasakian geometry
of links while reviewing basic facts about their differential geometry and topology.
At the end of this section we are left with a powerful
method of producing positive Sasakian structure on links.

In Section 7 we begin to discuss the famous Calabi Conjecture proved
in 1978 by Yau \cite{Yau78}. Yau's proof uses the continuity method and
works equally well for compact orbifolds. This fact has
important consequences for Sasakian manifolds: every positive
Sasakian manifold admits a metric of positive Ricci curvature. This observation
offers a very effective tool of proving the existence of such metrics on
many odd-dimensional manifolds. One interesting example is a theorem
of Wraith \cite{Wra97}, proved originally by surgery methods, which asserts that all homotopy
spheres which bound parallelizable manifolds admit metrics of positive Ricci curvature. The
authors together with M. Nakamaye \cite{BGN03b} recently gave an independent proof of this
result using the methods described here.

In the next section we turn our attention to positive K\"ahler-Einstein metrics. Even
in the smooth Fano category tractable necessary and sufficient
conditions for such a metric to exist are not known. After disproving one of the
Calabi conjectures asserting that in the absence of holomorphic vector fields
K\"ahler-Einstein metrics should exist, Tian proposed his own conjecture
\cite{Tia97} proving it in one direction.
Even assuming the conjecture to be true, in general it is
not easy to check if a particular Fano manifold (orbifold) satisfies the
required stability conditions.
On the other hand, in some cases, the continuity method has been used
effectively to check sufficient conditions. In this respect the method
of Demailly and Koll\'ar mentioned before draws
on earlier results of Nadel \cite{Nad90}, \cite{Siu88}, Tian \cite{Tia87, Tia89, Tia90}, and Tian
and Yau
\cite{TiY87}.

The last two sections give a summary of what has been accomplished to date by
applying the continuity method to Fano orbifolds.
We begin with a brief discussion of the method itself.
We follow with two important examples of how the method
applies to orbifolds constructed as hypersurfaces in weighted
projective spaces. We review our recent work \cite{BGK03,BGKT03} in collaboration with 
J. Koll\'ar which shows that standard odd-dimensional spheres have Einstein metrics
with one-dimensional isometry group and very large moduli spaces.
There we also proved that all homotopy spheres in dimensions $4n+1,7,11,$ and $15$ that 
bound 
parallelizable manifolds admit Sasakian-Einstein metrics. 
We discuss a conjecture that the last statement is true in
any odd dimension. Furthermore, it is shown that in each odd dimension starting with
$n=5$ there are infinitely many rational homology spheres which admit
Einstein metrics \cite{BG03p}. We close with the discussion of Sasakian-Einstein
geometry of Barden manifolds.

\begin{ack} We would like to than J\'anos Koll\'ar for comments. The authors were
partially supported by the NSF under grant number DMS-0203219. KG
would also like to thank Max-Planck-Institut f\"ur Mathematik in
Bonn for hospitality and support. This paper was written during
his one year visit there. In addition, KG would like to thank the
organizers of the XXIVth Workshop on {\it Geometry and Physics}
for support and Erwin Schr\"odinger Institute in Vienna for
hospitality and support during his short visit there.
\end{ack}

\bigskip
\section{Contact Structures and Sasakian Metrics}
\bigskip
Contact transformations arose in the theory of Analytical Mechanics developed in
the 19th century by Hamilton, Jacobi, Lagrange, and Legendre. But its first systematic
treatment was given by Sophus Lie. Consider $\bbr^{2n+1}$ with Cartesian coordinates
$(x^1,\cdots,x^n;y^1,\cdots,y^n;z),$ and a 1-form $\eta$ given by
\begin{equation}\label{contacteqn}
\eta=dz-\sum_iy^idx^i
\end{equation}
It is easy to see that $\eta$ satisfies $\eta\wedge (d\eta)^n\neq 0$.
A 1-form on $\bbr^{2n+1}$ that satisfies this equation
is called a {\it contact form}. Locally we have the following
\medskip\noindent

\begin{thm}  Let $\eta$ be a 1-form on $\bbr^{2n+1}$
that satisfies $\eta\wedge (d\eta)^n\neq 0$.
Then there is an open set $U\subset \bbr^{2n+1}$ and local coordinates
$(x^1,\cdots,x^n;y^1,\cdots,y^n;z)$ such that $\eta$ has the form (\ref{contacteqn}) in $U.$
\end{thm}

\begin{defn} A $(2n+1)$-dimensional manifold $M$ is a
{\it contact manifold} if there
exists a 1-form $\eta$, called a {\it contact 1-form}, on $M$ such that
$$\eta\wedge (d\eta)^n \neq 0$$
everywhere on $M.$ A {\it contact structure} on $M$ is an equivalence class of
such
1-forms, where $\eta'\sim \eta$ if there is a nowhere vanishing function $f$ on
$M$ such
that $\eta'=f\eta.$
\end{defn}

\begin{lem}
On a contact manifold $(M,\eta)$ there is a unique vector field $\xi$,
called the {\it Reeb
vector field}, satisfying the two conditions
$$\xi \rfloor \eta=1, \qquad \xi \rfloor d\eta =0.$$
\end{lem}

\begin{defn}
An {\it almost contact
structure} on a differentiable manifolds $M$ is a triple $(\xi,\eta,\Phi),$ where $\Phi$ is a
tensor field of type $(1,1)$ (i.e. an endomorphism of $TM$), $\xi$ is a vector field, and
$\eta$
is a 1-form which satisfy
$$\eta(\xi)=1 ~~~\hbox{and}~~~ \Phi\circ \Phi= -\bbi + \xi\otimes\eta,$$
where $\bbi$ is the identity endomorphism on $TM.$ A smooth manifold with such a
structure is called an {\it almost contact manifold}.
\end{defn}

Let $(M,\eta)$ be a contact manifold with a contact 1-form $\eta$ and consider
$\cald = \ker~\eta\subset TM.$  The subbundle $\cald$ is maximally
non-integrable and it is called the {\it contact
distribution}. The pair $(\cald,\gro)$, where $\gro$ is the
restriction of $d\eta$ to $\cald$ gives $\cald$ the structure of a symplectic
vector bundle. We denote by $\calj(\cald)$ the space of all almost complex
structures $J$ on $\cald$ that are compatible with $\gro,$ that is the subspace
of smooth sections $J$ of the endomorphism bundle ${\rm End}(\cald)$
that satisfy
\begin{equation}
J^2= -\bbi, \ \ d\eta(JX,JY)=d\eta(X,Y),\ \
d\eta(X,JX)>0
\end{equation}
for any smooth sections $X,Y$ of $\cald.$ Notice
that each $J\in \calj(\cald)$ defines a Riemannian metric $g_\cald$ on $\cald$
by setting
\begin{equation}\label{tranmetric}
g_\cald(X,Y) =d\eta(X,JY).
\end{equation}
One easily checks that $g_\cald$
satisfies the compatibility condition $g_\cald(JX,JY)=g_\cald(X,Y).$
Furthermore, the map $J\mapsto g_\cald$ is one-to-one, and the space
$\calj(\cald)$ is contractible. A choice of $J$ gives $M$ an almost CR structure.

Moreover, by extending $J$ to all of $TM$ one obtains an almost contact
structure. There are some choices of conventions to make here. We
define the section $\Phi$ of ${\rm End}(TM)$ by $\Phi =J$ on $\cald$ and
$\Phi\xi=0$, where $\xi$ is the Reeb vector field associated to $\eta.$ We can
also extend the transverse metric $g_\cald$ to a metric $g$ on all of $M$ by
\begin{equation}\label{sametric}
g(X,Y)= g_\cald +\eta(X)\eta(Y)= d\eta(X,\Phi Y)+
\eta(X)\eta(Y)
\end{equation}
for all vector fields $X,Y$ on $M.$ One
easily sees that $g$ satisfies the compatibility condition $g(\Phi X,\Phi
Y)=g(X,Y)-\eta(X)\eta(Y).$

\begin{defn}
A contact manifold $M$ with a contact form
$\eta$, a vector field $\xi,$ a section
$\Phi$ of ${\rm End}(TM),$ and
a Riemannian metric $g$ which satisfy the conditions
$$\eta(\xi)=1,\qquad \Phi^2=-\bbi +\xi\otimes \eta,$$
$$g(\Phi X,\Phi Y)
=g(X,Y)-\eta(X)\eta(Y)$$
is known as a {\it metric contact structure} on $M.$
\end{defn}

\begin{defn-thm} A Riemannian manifold $(M,g)$ is called a
{\it Sasakian manifold} if any one, hence all, of the following equivalent
conditions hold:
\begin{enumerate}
\item
There exists a Killing vector field $\xi$ of unit length on $M$
so that the tensor field $\Phi$ of type $(1,1)$, defined by
$\Phi(X) ~=~ -\nabla_X \xi$,
satisfies the condition
$$(\nabla_X \Phi)(Y) ~=~ g(X,Y)\xi -g(\xi,Y)X$$
for any pair of vector fields $X$ and $Y$ on $M.$
\item There exists a Killing vector field $\xi$ of unit length on $M$
so that the Riemann curvature satisfies the condition
$$R(X,\xi)Y ~=~ g(\xi,Y)X-g(X,Y)\xi,$$
for any pair of vector fields $X$ and $Y$ on $M.$
\item The metric cone
$(\calc(M), \bar{g}) ~=~ (\bbr_+\times M, \ dr^2+r^2g)$ is K\"ahler.
\end{enumerate}
\end{defn-thm}

We refer to the quadruple $\cals=(\xi,\eta,\Phi,g)$ as a
{\it Sasakian structure} on $M$, where $\eta$ is the 1-form dual
vector field $\xi.$ It is easy to see that
$\eta$ is a contact form whose Reeb vector field is $\xi$.
In particular  $\cals=(\xi,\eta,\Phi,g)$ is a special type of
{\it metric contact structure}.

The vector field $\xi$ is nowhere vanishing, so
there is a 1-dimensional foliation $\calf_\xi$ associated with every Sasakian
structure, called the {\it characteristic foliation}. We will denote the
space of leaves of this foliation by $\calz$. Each leaf of  $\calf_\xi$
has a holonomy group associated to it.  The dimension of the
closure of the leaves is called the {\it rank} of $\cals$.  We
shall be interested in the case $\hbox{rk}(\cals)=1$. We have

\begin{defn}\label{quasi-regular}
The characteristic foliation $\calf_\xi$ is said to be {\it
quasi-regular}\index{quasi-regular} if there is a positive integer
$k$ such that each point has a foliated coordinate chart $(U,x)$
such that each leaf of $\calf_\xi$ passes through $U$ at most $k$
times. Otherwise $\calf_\xi$ is called irregular. If $k=1$ then the foliation is called {\it
regular}\index{regular}, and we use the terminology {\it non-regular} to mean quasi-regular, but 
not regular.
\end{defn}

\bigskip
\section{Transverse K\"ahler Geometry}
\bigskip

Let $(M,\xi,\eta,\Phi,g)$ be a Sasakian manifold, and consider the
contact subbundle $\cald=\ker~\eta.$ There is an orthogonal
splitting of the tangent bundle as
\begin{equation}
TM=\cald \oplus L_\xi,
\end{equation}
where $L_\xi$ is the trivial line bundle generated by the Reeb
vector field $\xi.$ The contact subbundle $\cald$ is just the
normal bundle to the characteristic foliation $\calf_\xi$
generated by $\xi.$ It is naturally endowed with both a complex
structure $J=\Phi|\cald$ and a symplectic structure $d\eta.$
Hence, $(\cald,J,d\eta)$ gives $M$ a {\it transverse K\"ahler}
structure with K\"ahler form $d\eta$ and metric $g_\cald$ defined
as in (\ref{tranmetric})
which is related to the Sasakian metric $g$ by
$g=g_\cald \oplus \eta\otimes \eta$ as in (\ref{sametric}).
We have the following fundamental structure theorems
\cite{BG00a}:

\begin{thm}
Let $(M,\xi,\eta,\Phi,g)$ be a compact quasi-regular Sasakian manifold of
dimension $2n+1$, and let $\calz$ denote the space of leaves of the
characteristic foliation. Then
the leaf space $\calz$ is a Hodge orbifold with
K\"ahler metric $h$ and K\"ahler form $\gro$ which defines an integral class
$[\gro]$ in $H^2_{orb}(\calz,\bbz)$ so that $\pi:(M,g) \ra{1.3}
(\calz,h)$ is an orbifold Riemannian submersion. The fibers of $\pi$ are
totally geodesic submanifolds of $M$ diffeomorphic to $S^1.$
\end{thm}

\begin{thm} Let $(\calz,h)$ be a Hodge orbifold.
Let $\pi:M\ra{1.3} \calz$ be the $S^1$ V-bundle whose first Chern
class is $[\gro],$ and let $\eta$ be a connection 1-form in $M$ whose
curvature is $2\pi^*\gro,$ then $M$ with the metric
$\pi^*h+\eta\otimes\eta$ is a Sasakian orbifold.  Furthermore, if all the local
uniformizing groups inject into the group of the bundle $S^1,$ the total space
$M$ is a smooth Sasakian manifold.
\end{thm}

\begin{rem} The structure theorems discussed above show that
there are two K\"ahler geometries naturally associated with
every Sasakian manifold and we get the following diagram

\begin{equation}\label{sas}
\begin{matrix}
\calc(M)&\hookleftarrow&M \\
&{}&\decdnar\pi \\
{}&{}&\calz \\
\end{matrix}
\end{equation}
\end{rem}

The orbifold cohomology groups $H_{orb}^p(\calz,\bbz)$ were defined by
Haefliger \cite{Hae84}. In analogy with the smooth case a {\it Hodge orbifold}
is then defined to be a compact K\"ahler orbifold whose K\"ahler class
lies in $H_{orb}^2(\calz,\bbz)$. Alternatively, we can develop the
concept of basic cohomology. This is useful in trying to extend
the notion of $\calz$ being Fano to the
orbifold situation. This can be done in several ways. Here we will
use the notion of basic Chern classes.
Recall \cite{Ton} that a smooth p-form $\gra$ on $M$ is called
{\it basic} if
\begin{equation}
\xi\rfloor \gra=0, \qquad \call_\xi\gra=0,
\end{equation}
and we let $\grL^p_B$ denote the sheaf of germs of basic p-forms
on $M,$ and by $\grO_B^p$ the set of global sections of $\grL^p_B$
on $M.$ The sheaf $\grL^p_B$ is a module over the ring,
$\grL^0_B,$ of germs of smooth basic functions on $M.$ We let
$C^\infty_B(M)=\grO^0_B$ denote global sections of $\grL^0_B,$
i.e. the ring of smooth basic functions on $M.$  Since exterior
differentiation preserves basic forms we get a de Rham complex
\begin{equation}
\cdots\ra{2.5}\grO_B^p\fract{d}{\ra{2.5}}\grO_B^{p+1}\ra{2.5}\cdots
\end{equation}
whose cohomology $H^*_B(\calf_\xi)$ is called the {\it basic
cohomology} of $(M,\calf_\xi).$ The basic cohomology ring
$H^*_B(\calf_\xi)$ is an invariant of the foliation $\calf_\xi$
and hence, of the Sasakian structure on $M.$ It is related to the
ordinary de Rham cohomology $H^*(M,\bbr)$ by the long exact
sequence \cite{Ton}
\begin{equation}\label{exact}
\cdots\ra{2.5}H_B^p(\calf_\xi)\ra{2.5}H^p(M,\bbr)\fract{j_p}{\ra{2.5}}
H_B^{p-1}(\calf_\xi) \fract{\grd}{\ra{2.5}}
H^{p+1}_B(\calf_\xi)\ra{2.5}\cdots
\end{equation}
where $\grd$ is the connecting homomorphism given by
$\grd[\gra]=[d\eta\wedge \gra]=[d\eta]\cup[\gra],$ and $j_p$ is
the composition of the map induced by $\xi\rfloor$ with the well
known isomorphism $H^r(M,\bbr)\approx H^r(M,\bbr)^{S^1}$ where
$H^r(M,\bbr)^{S^1}$ is the $S^1$-invariant cohomology defined from
the  $S^1$-invariant r-forms $\grO^r(M)^{S^1}.$ We also note that
$d\eta$ is basic even though $\eta$ is not.
Next we exploit the fact that the transverse geometry is K\"ahler.
Let $\cald_\bbc$ denote the complexification of $\cald,$ and
decompose it into its eigenspaces with respect to $J,$ that is,
$\cald_\bbc= \cald^{1,0}\oplus \cald^{0,1}.$ Similarly, we get a
splitting of the complexification of the sheaf $\grL^1_B$ of basic
one forms on $M,$ namely
$$\grL^1_B\otimes \bbc = \grL^{1,0}_B\oplus \grL^{0,1}_B.$$
We let $\cale^{p,q}_B$ denote the sheaf of germs of basic forms of
type $(p,q),$ and we obtain a splitting
\begin{equation}
\grL^r_B\otimes \bbc = \bigoplus_{p+q=r}\cale^{p,q}_B.
\end{equation}

The basic cohomology groups $H^{p,q}_B(\calf_\xi)$ are fundamental
invariants of a Sasakian structure which enjoy many of the same
properties as the ordinary Dolbeault cohomology of a K\"ahler
structure.

Consider the complex vector bundle $\cald$ on a Sasakian manifold
$(M,\xi,\eta,\Phi,g).$ As such $\cald$ has Chern classes
$c_1(\cald),\cdots,c_n(\cald)$ which can be computed by choosing a
connection $\nabla^\cald$ in $\cald$ \cite{Kobbook}. Let us choose
a local  foliate unitary transverse frame $(X_1,\cdots,X_n),$ and
denote by $\grO^T$ the transverse curvature 2-form with respect to
this frame. A simple calculation shows that
$\grO^T$ is a basic $(1,1)$-form.
Since the curvature 2-form $\grO^T$ has type $(1,1)$ it follows as
in ordinary Chern-Weil theory that

\begin{defn-thm}\label{transChernWeil}
The $k$th Chern class $c_k(\cald)$ of the complex vector bundle
$\cald$ is represented by the basic $(k,k)$-form $\grg_k$
determined by the formula
$$\det(\bbi_n-\frac{1}{2\pi i}\grO^T)=1+\grg_1+\cdots +\grg_k.$$
Since $\grg_k$ is a closed basic $(k,k)$-form it represents an element in
$H_B^{k,k}(\calf_\xi)\subset H_B^{2k}(\calf_\xi)$ that is called
the {\it basic $k$th Chern class} and denoted by $c_k(\calf_\xi).$
\end{defn-thm}

We now concentrate on the first Chern classes $c_1(\cald)$ and
$c_1(\calf_\xi).$ We have

\begin{defn}\label{c1def}
A Sasakian structure $(\xi,\eta,\Phi,g)$ is said to be {\it
positive (negative)} if $c_1(\calf_\xi)$ is represented by a
positive (negative) definite $(1,1)$-form. If either of these two
conditions is satisfied $(\xi,\eta,\Phi,g)$ is said to be {\it
definite}, and otherwise $(\xi,\eta,\Phi,g)$ is called {\it
indefinite}. $(\xi,\eta,\Phi,g)$ is said to be {\it null} if
$c_1(\calf_\xi)=0.$
\end{defn}

In analogy with common terminology of smooth algebraic varieties
we see that a positive Sasakian structure is a {\it transverse
Fano structure}, while a null Sasakian structure is a {\it
transverse Calabi-Yau structure}. The negative Sasakian case
corresponds to the canonical bundle being ample; we refer to this
as a {\it transverse canonical structure}.

\begin{rem}\label{Fano1} Alternatively, a complex orbifold $\calz$ is Fano if its
orbifold canonical bundle $K_{\calz^{orb}}$ is anti-ample. In the
case $\calz$ is well-formed, that is when the orbifold singularities
have codimension at least 2, the orbifold canonical bundle $K_{\calz^{orb}}$ can be identified 
with 
the ordinary canonical bundle.
However, in the presence of codimension 1 singularities the orbifold
canonical divisor {\bf is not} the usual algebraic geometric
canonical divisor, but is shifted off by the ramification divisors coming from the codimension one 
singularities \cite{BGK03}. We shall give specific examples of this difference
later.
\end{rem}

\bigskip
\section{The Einstein Condition}
\bigskip
\begin{defn}
A Sasakian space $(M,g)$ is {\it Sasakian-Einstein} if the metric $g$ is also
Einstein. For any 2n+1-dimensional Sasakian manifold
${\rm Ric}(X,\xi)=2n\eta(X)$ implying that any Sasakian-Einstein
metric must have positive scalar
curvature. Thus any complete Sasakian-Einstein
manifold must have a finite
fundamental group. Furthermore
the metric cone on $M$
$(\calc(M),\bar{g})= (\bbr_+\times M, \ dr^2+r^2g)$ is K\"ahler Ricci-flat (Calabi-Yau).
\end{defn}

The following theorem \cite{BG00a} is an
orbifold version of the famous Kobayashi bundle construction
of Einstein metrics on bundles over positive
K\"ahler-Einstein manifolds \cite{Bes, Kob56}.

\begin{thm}\label{kob}
 Let $(\calz, h)$ be a compact Fano orbifold
with $\pi_1^{orb}(\calz)=0$ and K\"ahler-Einstein metric $h$.
Let $\pi:M\ra{1.3} \calz$ be the $S^1$
V-bundle whose first Chern class is $\frac{c_1(\calz)}{\hbox{Ind}(\calz)}.$
Suppose further that the local uniformizing groups of $\calz$ inject into
$S^1.$ Then with the metric $g=\pi^*h+\eta\otimes\eta$,
$M$ is a compact simply connected Sasakian-Einstein manifold.
\end{thm}

Here ${\hbox{Ind}(\calz)}$ is the {\it orbifold Fano index} \cite{BG00a} defined to be the largest 
positive integer such that $\frac{c_1(\calz)}{\hbox{Ind}(\calz)}$ defines a class in the orbifold 
cohomology group $H^2_{orb}(\calz,\bbz).$  A very special class of Sasakian-Einstein spaces is 
naturally
related to several quaternionic geometries.

\begin{defn}
Let $(M,g)$ be a Riemannian manifold of dimension $m$.
We say that
$(M,g)$ is 3-Sasakian if the
metric cone $(\calc(M),\bar{g})=
(\bbr_+\times\cals, \ dr^2+r^2g)$ on $M$  is hyperk\"ahler.
\end{defn}

\begin{rem} In the 3-Sasakian case there is an extra structure, i.e.,
the transverse geometry $\calo$ of
the 3-dimensional foliation which is quaternionic-K\"ahler. In this case,
the transverse space $\calz$ is the twistor space of $\calo$ and
the natural map $\calz\ \ra{1.3}\ \calo$ is the orbifold twistor fibration \cite{Sal82}.
We get the following diagram which we denote by $\Diamond(M)$ and which extends
the diagram in (\ref{sas}) \cite{BGM93a, BGM94a}:
\vskip 20pt
\begin{center}
\begin{tabular}{|c|}\hline {\bf Hyperk\"ahler}\\{\bf Geometry}\\ \hline\end{tabular}
\end{center}
\begin{equation}
\begin{tabular}{|c|}\hline {\bf Twistor}\\{\bf Geometry}\\ \hline
\end{tabular}
\begin{array}{ccccc}
&&\calc(M)&&\\
&\swarrow&&\searrow&\\
\calz &&\hskip -15pt\la{4}\hskip -30pt\decdnar{}&& M\\
&\searrow& &\swarrow&\\
&&\calo&&
\end{array}
\begin{tabular}{|c|}\hline {\bf 3-Sasakian}\\{\bf Geometry}\\ \hline
\end{tabular}
\end{equation}
\begin{center}
\begin{tabular}{|c|}\hline {\bf Quaternion K\"ahler}\\{\bf Geometry}\\ \hline\end{tabular}
\end{center}
\end{rem}

\begin{rem} The table below summarizes properties of cone and transverse
geometries associated to various metric contact structures.
\vskip 12pt
\begin{center}
  \begin{tabular}{|c|c|c|}
  \hline
    Cone Geometry of $\calc(M)$ & $M$ & Transverse Geometry of $\calf_\xi$ \\ \hline\hline
    Symplectic & Contact & Symplectic \\ \hline
    K\"ahler &Sasakian & K\"ahler \\ \hline
    K\"ahler &positive Sasakian & Fano, $c_1(\calz)>0$ \\ \hline
    K\"ahler &null Sasakian & Calabi-Yau, $c_1(\calz)=0$ \\ \hline
    K\"ahler &negative Sasakian & canonical, $c_1(\calz)<0$ \\ \hline
    Calabi-Yau &Sasakian-Einstein & Fano, K\"ahler-Einstein \\ \hline
    Hyperk\"ahler &3-Sasakian& $\bbc$-contact, Fano, K\"ahler-Einstein \\ \hline
  \end{tabular}
\end{center}
\end{rem}

\bigskip
\section{Some Examples}
\bigskip
Below we list some well-know constructions of Sasakian-Einstein and 3-Sasakian
manifolds. We start with the latter.

\begin{exmp} Examples of 3-Sasakian manifolds are numerous and they are easily
constructed by way of the so called 3-Sasakian reduction \cite{BGM94a}. To begin with one
starts
with the canonical example of $\Diamond(M)$ where $M$ is the round $(4n-1)$-dimensional
sphere
$S^{4n-1}$ of constant sectional curvature 1.
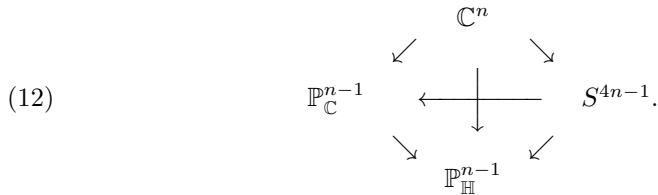
\begin{equation}\label{3sflat}
\begin{array}{ccccc}
&&\bbc^n&&\\
&\swarrow&&\searrow&\\
\bbp_\bbc^{n-1} &&\hskip -15pt\la{4}\hskip -30pt\decdnar{}&& S^{4n-1}.\\
&\searrow& &\swarrow&\\
&&\bbp_\bbh^{n-1}&&
\end{array}
\end{equation}

Moreover, there is such a diamond diagram for any semisimple Lie group $G$ and we get
all homogeneous examples this way \cite{BGM94a}, i.e.,

\begin{equation}\label{3shom}
\begin{array}{ccccc}
&&\calc(G/L)&&\\
&\swarrow&&\searrow&\\
\frac{G}{U(1)\cdot L}&&\hskip -15pt\la{4}\hskip -30pt\decdnar{}&& G/L.\\
&\searrow& &\swarrow&\\
&&\frac{G}{Sp(1)\cdot L}&&
\end{array}
\end{equation}
The spaces $\frac{G}{Sp(1)\cdot L}$ are the well-known Wolf spaces \cite{Wol65} and they are
the only know examples of smooth compact positive quaternion K\"ahler manifolds.
A conjecture of LeBrun and Salamon \cite{LeSal} asserts that there are no other examples.
This conjecture have been proved in the first three quaternionic dimensions by
Hitchin \cite{Hit74a}, Friedrich and Kurke \cite{FrKur}, Poon and Salamon
\cite{PoSal91}, and Herrera and Herrera \cite{HeHe02a, HeHe02b}.
\end{exmp}

\begin{exmp}\label{3stor} Now, one can start with any of the homogeneous diamonds
$\Diamond(G/L)$ in \ref{3shom}. In principle, any subgroup $H\subset G$ of appropriate
dimension leads to a reduction of $\Diamond(G/L)$ by symmetries of $H$. At
various levels of the diamond such reductions are known as hyperk\"ahler \cite{HKLR},
3-Sasakian
\cite{BGM94a}, and
quaternionic K\"ahler quotients \cite{GL}, respectively. In practice, it is not
easy to assure that, say, the 3-Sasakian quotient of $G/L$ by $H$ be a smooth
manifold. On the other hand, there are many cases when this happens.
For instance,
one can reduce the standard diagram \ref{3sflat} by an action of
$T^k\subset T^n\subset Sp(n)$
$k$-dimensional torus. We get the reduction diagram
\begin{equation*}
\begin{array}{ccccc}
&&\bbc^n&&\\
&\swarrow&&\searrow&\\
\bbp_\bbc^{n-1} &&\hskip -15pt\la{4}\hskip -30pt\decdnar{}&& S^{4n-1}\\
&\searrow& &\swarrow&\\
&&\bbp_\bbh^{n-1}&&
\end{array}
\begin{tabular}{|c|}\hline {\bf Reduction}\\ $\Longrightarrow$ \\ {\bf by k-torus}\\ \hline
\end{tabular}
\begin{array}{ccccc}
&&\calc(\cals(\Omega))&&\\
&\swarrow&&\searrow&\\
\calz(\Omega)&&\hskip -15pt\la{4}\hskip -30pt\decdnar{}&& \cals(\Omega).\\
&\searrow& &\swarrow&\\
&&\calo(\Omega)&&
\end{array}
\end{equation*}

The geometry of the 3-Sasakian reduced space $\cals(\Omega)$
is completely determined
by an integral $k\times n$ matrix $\Omega$ which for any $k\geq1$ defines a homomorphism
 $f_\Omega:T^k\ra{1.2} T^n\subset U(n)\subset Sp(n)$. The real dimension
${\rm dim}(\cals(\Omega))=4(n-k)-1$. In the  case ${\rm dim}(\cals(\Omega))=7$, there
are choices of $\Omega$ for any $k\geq1$ which make $\cals(\Omega)$ smooth. Since
$b_2(\cals(\Omega))=k$ we conclude that in dimension 7 there exist Einstein manifolds
with arbitrarily large second Betti number. These were the first such examples and they
were constructed in \cite{BGMR}. Interestingly, this toric reduction does not give smooth 
manifolds 
with large second Betti numbers in dimensions greater than 7 \cite{BGM98a}. Nevertheless, one 
can obtain Sasakian-Einstein manifolds with arbitrary second Betti number in any odd 
dimension 
greater than seven by the join construction discussed in 31 below. Later Bielawski showed that 
in 
all allowed dimensions all toric examples
must occur through the above procedure \cite{Bi3}.  
\end{exmp}

\begin{exmp} The first non-toric examples in dimension 11 and 15
were obtained by Boyer, Galicki, and Piccinni \cite{BGP02}. These are toric quotients
of the diamond digram $\Diamond(G/L)$ for $G=SO(n), L=SO(4)\times SO(n-4)$.
Alternatively, these
can be thought of as non-Abelian reductions of $\Diamond(S^{4n-1})$.
\end{exmp}

\begin{exmp}
Recently, the first non-toric examples in dimension 7 were obtained by
Grove, Wilking, and Ziller \cite{zilpriv}. They use an
orbifold bundle construction with the examples of orbifold twistor space and self-dual Einstein
metrics $\calz_k\ra{1.2}\calo_k$
discovered by Hitchin in 1992 \cite{Hit92}.
The self-dual Einstein metric on $\calo_k$ is defined on
$S^4\setminus\bbr\bbp^2$
and it has $\bbz_k$ orbifold singularity along $\bbr\bbp^2$. However, it turns out that the
bundle $M_k\ra{1.2}\calz_k$ is actually smooth.
In particular, one can compute the integral cohomology ring of $M_k$. For odd $k$ the 
3-Sasakian manifold
$M_k$ is a  rational homology 7-sphere with non-zero torsion depending on $k$.
Hence, there exist infinitely many rational homology 7-spheres which have 3-Sasakian metrics
\end{exmp}

Let us turn our attention to complete examples of Sasakian-Einstein manifolds which are not
3-Sasakian.
\begin{exmp}
The standard example is that of complex Hopf fibration
\begin{equation}
\bbc^{n}\hookleftarrow S^{2n-1}\ \ra{1.2}\ \bbp_\bbc^{n-1}.
\end{equation}
Just as in 3-Sasakian case this example generalizes when one replaces the
complex projective space with a generalized flag manifold. That is, consider
any complex semi-simple Lie group $G$. A maximal solvable subgroup $B$ of $G$ is called a
{\it Borel subgroup} and is unique up to conjugacy. Any $P\subset G$ containing $B$ is
called {\it parabolic}. It is known that any such generalized flag $G/P$
admits a homogeneous K\"ahler-Einstein metric, and that any compact homogeneous simply
connected K\"ahler manifold is a generalized flag manifold. Applying the construction
of Theorem \ref{kob} gives all compact homogeneous Sasakian-Einstein metrics, in fact, all 
compact homogeneous Sasakian manifolds \cite{BG00a}.
\end{exmp}

\begin{exmp} The Kobayashi bundle construction also gives  many inhomogeneous examples.
These are all circle bundles over compact smooth Fano manifolds. For instance,
in the case of surfaces all del Pezzo surfaces are classified and it is known
which of them admit K\"ahler-Einstein metrics \cite{Tia90, Tia99, Tia00}.
\end{exmp}

\begin{thm}
The following del Pezzo surfaces admit
K\"ahler-Einstein metrics: $\bbc\bbp^2, \bbc\bbp^1\times\bbc\bbp^1,
\bbc\bbp^2\#n\overline{\bbc\bbp^2}, \ \ 3\leq n\leq8$.
Furthermore, the moduli space of K-E structures
in each case is completely understood.
\end{thm}

As and immediate consequence, for 5-manifolds, we have the following
result of Friedrich and Kath \cite{FrKat1}

\begin{thm} Let $M_l=S^5\#l(S^2\times S^3).$
\begin{enumerate}
\item For each $l=0,1,3,4,$ there is precisely one
regular Sasakian-Einstein
structure
on $M_l.$
\item For each $5\leq l\leq 8$ there is a $2(l-4)$ complex parameter
family of
inequivalent regular Sasakian-Einstein structures on $M_l.$
\item For $l=2$ or $l\geq 9$ there are no regular Sasakian-Einstein
structures on
$M_l.$
\end{enumerate}
\end{thm}

There are two del Pezzo surfaces which do not admit any K-E metrics due to
theorem of Matsushima \cite{Mat57}: the existence is obstructed by holomorphic
vector fields. These are blow-ups of $\bbc\bbp^2$ at one or two points.

\begin{rem} A well known result of
Martinet says that every orientable 3-manifold admits a contact structure.
Furthermore, all Sasakian 3-manifolds have been classified \cite{Bel00,Bel01,Gei97} and
they are Seifert bundles over Riemann
surfaces. In this case every compact Sasakian 3-manifold is
either negative, null, or positive. In addition, if $M$ is Sasakian-Einstein
than it follows that the universal cover $\tilde{M}$ is isomorphic to
the standard Sasakian-Einstein metric on $S^3$.
\end{rem}

\begin{exmp} \label{BarMan} [{\bf Barden Manifolds}]  Similarly one might try to classify all
Sasakian manifolds in dimension 5. In the simply connected case there exists a classification
result of all smooth 5-manifolds due to Smale \cite{Sm62} and  Barden \cite{Bar65}.
Extending Smale's theorem for spin manifolds Barden proves  the following:
\begin{thm}\label{Barden} The class of simply connected, closed, oriented,
smooth, 5-manifolds is classifiable under diffeomorphism. Furthermore,
any such $M$ is diffeomorphic to one of the spaces
$M_{j;k_1,\ldots,k_s}=X_j\#M_{k_1}\#\cdots \#M_{k_s}$,
where $-1\leq j\leq\infty,$
$s\geq0$, $1<k_1$ and $k_i$ divides $k_{i+1}$ or $k_{i+1}=\infty$.
A complete set
of invariants is provided by $H_2(M,\bbz)$ and an additional
diffeomorphism invariant $i(M)=j$ which depends only on the
second Stiefel-Whitney class $w^2(M)$.
\end{thm}
In this article we will refer to a simply connected, closed, oriented,
smooth, 5-manifold as a {\it Barden manifold}. The building blocks of
Theorem \ref{Barden}  are given in the table below.
They are listed with $H_2(M,\bbz)$ and Barden's $i(M)$ invariant.

\begin{center}
\begin{tabular}{|c|c|c|}\hline
$M$ & $H_2(M,\bbz)$ &$i(M)$ \\
\hline\hline
$X_{-1}=SU(3)/SO(3)$ & $\bbz_2$ &$-1$ \\ \hline
$X_{n},\ n\geq1$ & $\bbz_{2^n}\oplus\bbz_{2^n}$ &$n$ \\ \hline
$X_{\infty}=$ non-trivial $S^3$ bundle over $ S^2$ & $\bbz$ &$\infty$ \\ \hline
$X_{0}=S^5$ & $0$ &$0$ \\ \hline
$M_\infty=S^3\times S^2$ & $\bbz$ &$0$ \\ \hline
$M_{n}, \ n>1$ & $\bbz_n\oplus\bbz_n$ &$0$ \\ \hline
\end{tabular}
\end{center}
\vskip 12pt
When $M$ is spin $i(M)=j=0$ as is
$w^2(M)=0$ and
Barden's result is the extension of the well-known theorem of Smale for
spin 5-manifolds.
By an old theorem of Gray \cite{Gra59} $M$ admits an almost contact structure
when $j=0,\infty$ and by another result of Geiges $M$
is in such a case necessarily contact \cite{Gei91}.
\end{exmp}

\begin{ques} It is natural to ask whether Barden manifolds admit Sasakian structures.
Specifically, we would like to ask
\begin{enumerate}
\item Does every Barden manifold which is contact $(j=0,\infty)$ admit a Sasakian structure?
\item Does every Barden manifold which is spin  $(j=0)$
admit positive (respectively negative)
Sasakian structure? (Null Sasakian structures are obstructed. For example, Corollary 1.10 of 
\cite{BGN03a} implies that $S^5$ and $S^2\times S^3$ cannot admit null Sasakian structures).
\item Which of the Barden manifolds which are spin admit Sasakian-Einstein structures?
\end{enumerate}
\end{ques}

\begin{rem} There is one more construction of non-regular Sasakian-Einstein
manifolds which draws on examples of 3-Sasakian structures.
In \cite{BG00a} the authors observed that the set of all Sasakian-Einstein manifolds
has a monoid structure, i.e., for any two compact quasi-regular
Sasakian-Einstein orbifolds $M_1$ and $M_2$ one can define $M_1\star M_2$
which is automatically a compact Sasakian-Einstein orbifold of dimension
${\rm dim}(M_1)+{\rm dim}(M_2)-1$. This construction is an orbifold
generalization of the well-known construction of Wang and Ziller in \cite{WaZi90}
adapted to the Sasakian-Einstein setting. It turns out that a join of any non-regular
Sasakian-Einstein manifold $M_1$ with a regular Sasakian-Einstein space $M_2$ (say, for 
example, $M_2=S^3$) is automatically a compact smooth, non-regular Sasakian-Einstein
manifold. The join construction produces new examples beginning in dimension 7. If, however, 
$M_1$ is 3-Sasakian, new examples begin in dimension 9.
In addition, a completely new construction of inhomogeneous Sasakian-Einstein metrics
have been considered by Gauntlett, Martelli, Sparks, and Waldram \cite{GMSW04a, GMSW04b}.
Their metrics are in fact very explicit and they were first obtained indirectly by
considering general supersymmetric solutions in certain $D=11$ supergravity theory.
\end{rem}
\bigskip
\section{Sasakian Geometry of Links}
\bigskip
Consider the affine space $\bbc^{n+1}$ together with a weighted
$\bbc^*$-action given by $(z_0,\ldots,z_n)\mapsto
(\grl^{w_0}z_0,\ldots,\grl^{w_n}z_n),$ where the {\it weights} $w_j$ are
positive integers. It is convenient to view the weights as the components of a
vector $\bfw\in (\bbz^+)^{n+1},$ and we shall assume that
$\gcd(w_0,\ldots,w_n)=1.$
\begin{defn} We say that $f$ is a weighted homogeneous polynomial with weights
$\bfw$ and of degree $d$ if
$f\in \bbc[z_0,\ldots,z_n]$ and satisfies
\begin{equation}
f(\grl^{w_0}z_0,\ldots,\grl^{w_n}z_n)=\grl^df(z_0,\ldots,z_n).
\end{equation}
\end{defn}
We shall assume that the
origin in $\bbc^{n+1}$ is an isolated singularity.
\begin{defn}\label{link} The link of $f$ is
defined by
\begin{equation}
L_f= \{f=0\}\cap S^{2n+1},
\end{equation}
where $S^{2n+1}$
is the unit sphere in $\bbc^{n+1}.$
\end{defn}

\begin{rem}
$L_f$ is endowed with a natural
quasi-regular Sasakian structure \cite{Tak,YK,BG01b} inherited as a Sasakian submanifold of 
the
sphere $S^{2n+1}$ with  its ``weighted'' Sasakian structure
$(\xi_\bfw,\eta_\bfw,\Phi_\bfw,g_\bfw)$ which in the standard coordinates
$\{z_j=x_j+iy_j\}_{j=0}^n$ on $\bbc^{n+1}=\bbr^{2n+2}$ is determined by
\begin{equation}
\eta_\bfw = \frac{\sum_{i=0}^n(x_idy_i-y_idx_i)}{\sum_{i=0}^n
w_i(x_i^2+y_i^2)}, \qquad \xi_\bfw
=\sum_{i=0}^nw_i(x_i\partial_{y_i}-y_i\partial_{x_i}),
\end{equation}
and the standard Sasakian structure $(\xi,\eta,\Phi,g)$ on $S^{2n+1}.$

The quotient of $S^{2n+1}$ by the ``weighted $S^1$-action'' generated by the
vector field $\xi_\bfw$ is the weighted projective space
$\bbp(\bfw)=\bbp(w_0,\ldots,w_n),$ and we
have a commutative diagram:
\begin{equation}
\begin{array}{cccc}
L_f &\ra{2.5}& S^{2n+1}_\bfw&\\
  \decdnar{\pi}&&\decdnar{} &\\
   \calz_f &\ra{2.5} &\bbp(\bfw),&
\end{array}
\end{equation}
where the horizontal arrows are Sasakian and K\"ahlerian embeddings,
respectively, and the vertical arrows are orbifold Riemannian submersions.
$L_f$ is the total space of the principal $S^1$ V-bundle over the orbifold
$\calz_f$. Alternatively we will sometimes denote $\calz_f$ as $X_d\subset\bbp(\bfw)$
to indicate the weights and the degree of $f$. In such case we will also
write $L_f=L(X_d\subset\bbp(\bfw))$.
\end{rem}

\begin{prop}\label{FanoLink} \cite{BGK03}
The orbifold $\calz_f$ is Fano if and only if $d-\sum w_i<0$.
\end{prop}

\begin{exmp}\label{Fano2} At this point we  return to the comments made in
Remark \ref{Fano1}. Consider links defined by
\begin{equation}
f_k(z_0,z_1,z_2)=z_0^{6k-1}+z_1^3+z_2^2.
\end{equation}
The orbifold  $\calz_{f_k}$ is a hypersurface $X_d$ in
$\bbp(6,2(6k-1),3(6k-1))$ of degree $d=6(6k-1)$. The corresponding
link $L_{f_k}$ is of Brieskorn-Pham type and will be denoted by
$L(6k-1,3,2)$ (see (\ref{BP})).  All 3-dimensional Brieskorn-Pham
links were classified by Milnor in \cite{Mil75}. According to the
Proposition \ref{FanoLink}, $L(6k-1,3,2)$ is positive only when
$k=1$ and in all other cases it is negative. Indeed,
$L(5,3,2)\simeq S^3/I^*$ is the famous Poincar\'e sphere, where
$I^*\subset SU(2)$ is the binary isocahedral group. For $k>1$, the
link $L(6k-1,3,2)$ is a homology sphere with infinite fundamental
group. The complex orbifold $\calz_{f_k}$, for $k>1$ is not Fano.
In particular, it cannot have an orbifold metric of {\bf constant
positive} curvature (though it has a natural metric of constant
negative curvature). On the other hand, as an algebraic variety,
for any $k$ we must have $\calz_{f_k}\simeq\bbp^1$. This can be
seen from the generalized genus formula. For any curve
$X_d\subset\bbp(w_0,w_1,w_2)$ we have:
\begin{equation}
g(X_d)=\frac{1}{2}\Bigl(\frac{d^2}{w_0w_1w_2}-d\sum_{i<j}\frac{\gcd(w_i,w_j)}{
w_iw_j}+\sum_i\frac{\gcd(d,w_i)}{w_i}-1\Bigr).
\end{equation}
Hence,  $\calz_{f_k}$ is certainly Fano as a smooth variety in the 
algebraic geometric sense, but it has codimension1 orbifold singularities and it
{\bf is not} Fano in the orbifold sense. Here the orbifold
canonical class is not the usual algebraic geometric canonical
class, but the codimension one orbifold ramification divisors are added in. By
Milnor's classification $L(6k-1,3,2),$ for $(k>1),$ is the
quotient of the universal cover $\widetilde{SL}(2,\bbr)$ of
$SL(2,\bbr)$ by a co-compact discrete subgroup $\grG\subset
\widetilde{SL}(2,\bbr).$ Furthermore, $L(6k-1,3,2)$ has a finite
covering by a manifold that is diffeomorphic to a non-trivial
circle bundle over a Riemann surface of some genus $g>1.$
\end{exmp}

Now, recall the well-known construction of Milnor for isolated hypersurface
singularities \cite{Mil68, MiOr70}: there is a fibration of $(S^{2n+1}-L_f)\ra{1.3} S^1$ whose fiber
$F$ is an open manifold that is homotopy equivalent to a bouquet of n-spheres
$S^n\vee S^n\cdots \vee S^n.$ The {\it Milnor number} $\mu$ of $L_f$ is the
number of $S^n$'s in the bouquet. It is an invariant of the link which
can be calculated explicitly in terms of the degree $d$ and weights
$(w_0,\ldots,w_n)$ by the formula
\begin{equation}
\mu =\mu(L_f)=\prod_{i=0}^n\bigl(\frac{d}{w_i}-1\bigr)
\end{equation}
The closure $\bar{F}$ of $F$ has the same homotopy type as $F$ and is a compact
manifold with boundary precisely the link $L_f.$ So the reduced homology of
$F$ and $\bar{F}$ is only non-zero in dimension $n$ and $H_n(F,\bbz)\approx
\bbz^\mu.$ Using the Wang sequence of the Milnor fibration together with
Alexander-Poincare duality gives the exact sequence
\begin{equation}\label{exact1}
0\!\!\rightarrow\!\! H_n(L_f,\bbz)\rightarrow H_n(F,\bbz)
\fract{\bbi -h_*}{\ra{1.5}} H_n(F,\bbz) \rightarrow H_{n-1}(L_f,\bbz)\!\!\rightarrow\!\!0
\end{equation}
where $h_*$ is the {\it monodromy} map (or characteristic
map) induced by the $S^1_\bfw$ action. From this we see that
$H_n(L_f,\bbz)=\ker(\bbi -h_*)$ is a free Abelian group, and $H_{n-1}(L_f,\bbz)
={\rm Coker}(\bbi -h_*)$ which in general has torsion, but whose free part equals
$\ker(\bbi -h_*).$ So the topology of $L_f$ is encoded in the monodromy map
$h_*.$  There is a well-known algorithm due to Milnor and Orlik \cite{MiOr70}
for computing the free part of $H_{n-1}(L_f,\bbz)$ in terms of the
characteristic polynomial $\grD(t)=\det(t\bbi -h_*),$ namely the Betti number
$b_n(L_f)=b_{n-1}(L_f)$ equals the number of factors of $(t-1)$ in $\grD(t).$
First we mention an important immediate consequence of the
exact sequence (\ref{exact1}) which is due to Milnor:
\medskip
\begin{prop}
The following hold:
\begin{enumerate}
\item $L_f$ is a rational homology sphere if and only if $\grD(1)\neq 0.$
\item $L_f$ is a homology sphere if and only if $|\grD(1)|=1.$
\item If $L_f$ is a rational homology sphere, then the order of $H_{n-1}(L_f,\bbz)$ equals
$|\grD(1)|.$
\end{enumerate}
\end{prop}
\begin{exmp} The following table lists some illustrating examples. All of
the chosen links are Fano, but negative and null Sasakian structures can
also be considered. We either explicitly identify the link with some
smooth contact manifold or list non-vanishing homology groups. $\Sigma_k^7$ and
$\Sigma_p^7$ indicate homotopy spheres where the differentiable structure depends on
$k$ and $p$.
\vskip 12pt
\begin{center}
\begin{tabular}{|c|c|}\hline
$\calz_f$ & $L_f$  \\
\hline\hline
$X_2\subset\bbp(1,1,1,1)$ & $S^2\times S^3$  \\ \hline
$X_3\subset\bbp(1,1,1,1)$ & $6\#(S^2\times S^3)$  \\ \hline
$X_4\subset\bbp(1,1,1,2)$ & $7\#(S^2\times S^3)$  \\ \hline
$X_6\subset\bbp(1,1,2,3)$ & $8\#(S^2\times S^3)$  \\ \hline
$X_{k+1}\subset\bbp(1,1,1,k)$ & $k\#(S^2\times S^3)$  \\ \hline
$X_{10}\subset\bbp(1,2,3,5)$ & $9\#(S^2\times S^3)$  \\ \hline
$X_{127}\subset\bbp(11,29,39,49)$ & $2\#(S^2\times S^3)$  \\ \hline
$X_{256}\subset\bbp(13,35,81,128)$ & $S^2\times S^3$  \\ \hline
$X_{3k}\subset\bbp(3,3,3,k)$, $k\not=3n$ & $M_k$  \\ \hline
$X_{4k}\subset\bbp(4,4,4,4,k)$, $k\not=2n$ & $|H_3(L_f,\bbz)|=k^{21}$  \\ \hline
$X_{6(6k-1)}\subset\bbp(6,2(6k-1), 3(6k-1),3(6k-1),3(6k-1))$ & $\Sigma_k\simeq S^7$  \\
\hline
$X_{2(2p+1)}\subset\bbp(2,2p+1, 2p+1,2p+1,2p+1,2p+1)$ & $\Sigma_p\simeq S^9$  \\ \hline
\end{tabular}
\end{center}
\vskip 12pt
\end{exmp}

In some of the above examples the homogeneous polynomial
$f$ can be chosen to contain no ``mixed" monomial terms of the form $z_i^{c_i}z_j^{c_j}$. Such
an $f$ is
called  of Brieskorn-Pham type. In his famous work, in 1966
Brieskorn considered links  $L(\bfa)$ defined by
\begin{equation}\label{BP}
\sum_{i=0}^{n}|z_i^2|=1,\qquad f_\bfa(\bfz)= z_0^{a_1}+\cdots +z_{n}^{a_n}=0.
\end{equation}

To the vector $\bfa=(a_0,\cdots,a_n)\in \bbz_+^{n+1}$ one associates a graph
$G(\bfa)$ whose $n+1$ vertices are labeled by $a_0,\cdots,a_n.$ Two vertices
$a_i$ and $a_j$ are connected if and only if $\gcd(a_i,a_j)>1.$ Let $C_{ev}$
denote the connected component of $G(\bfa)$ determined by the even integers.
Note that all even vertices belong to $C_{ev},$ but $C_{ev}$ may contain odd
vertices as well. Then we have the so-called {\bf Brieskorn Graph Theorem}
\cite{Bri66}:

\begin{thm}
The following hold:
\begin{enumerate}
\item The link $L(\bfa)$ is a rational homology sphere if and only if either
$G(\bfa)$ contains at least one isolated point, or $C_{ev}$ has an odd number
of vertices and for any distinct $a_i,a_j\in C_{ev},~$ $\gcd(a_i,a_j)=2.$
\item The link $L(\bfa)$ is an integral homology sphere if and only if either
$G(\bfa)$ contains at least two isolated points, or
$G(\bfa)$ contains one isolated point and $C_{ev}$ has an odd number
of vertices and $a_i,a_j\in C_{ev},~$ implies $\gcd(a_i,a_j)=2$ for any distinct  $i,j.$
\end{enumerate}
\end{thm}

Recall that by the seminal work of Milnor \cite{Mil56}, Kervaire and Milnor \cite{KeMil63}, and 
Smale
\cite{Sm61},
for each $n\geq 5$, differentiable
homotopy spheres of dimension $n$ form an Abelian group $\Theta_n$,
where the group operation is connected sum.
$\Theta_n$ has a subgroup $bP_{n+1}$ consisting of
those homotopy $n$-spheres which bound parallelizable manifolds $V_{n+1}.$
Kervaire and Milnor proved that $bP_{2m+1}=0$  for $m\geq 1,$
$bP_{4m+2}=0,$ or $\bbz_2$ and is $\bbz_2$ if $4m+2\neq 2^i-2$
for any $i\geq 3.$ The most interesting groups are
 $bP_{4m}$ for $m\geq 2.$ These are cyclic of order
\begin{equation}
|bP_{4m}|=2^{2m-2}(2^{2m-1}-1)~\hbox{numerator}~\biggl(\frac{4B_m}{m}\biggr),
\end{equation}
where $B_m$ is the $m$-th Bernoulli number.
Thus, for example $|bP_8|=28, |bP_{12}|=992,
|bP_{16}|=8128$ and $|bP_{20}|=130,816$. In the first two cases these
include all exotic spheres.
The correspondence is given by
\begin{equation}\label{exot}
KM: \Sigma\mapsto \frac{1}{8}\tau(V_{4m}(\Sigma)) {\rm mod} |bP_{4m}|,
\end{equation}
where $V_{4m}(\Sigma)$ is any parallelizable manifold bounding $\Sigma$
and $\tau$ is its signature.
Let $\Sigma_i$ denote the exotic sphere with $KM(\Sigma_i)=i$.
Now, the Brieskorn Graph Theorem tells us for which $\bfa$ the Brieskorn-Pham link
$L(\bfa)$ is a homotopy sphere. By (\ref{exot}) we need to be able
to compute the signature to determine  the diffeomorphism types of various links. We restrict our
interest just to the case when  $m=2k+1$.

In this case, the diffeomorphism type of a homotopy sphere
$L(\bfa)\in
bP_{2m-2}$ is determined by the signature $\grt(M)$ of a parallelizable manifold
$M$
whose boundary is $\Sigma_\bfa^{2m-3}.$ By the Milnor Fibration Theorem we can
take
$M$ to be the Milnor fiber $M_\bfa^{2m-2}$ which, for links of isolated singularities coming from
weighted homogeneous polynomials is diffeomorphic to the hypersurface  $\{\bfz\in
\bbc^{m} ~|~f_\bfa(z_0,\cdots,z_{m-1})=1\}.$

Brieskorn shows that the signature of $M^{4k}(\bfa)$ can be
written combinatorially as
\begin{eqnarray}\label{briescomb}
\notag \grt(M^{4k}(\bfa))&=& \#\bigl\{\bfx\in \bbz^{2k+1}
~|~0<x_i<a_i~\hbox{and}~0<\sum_{j=0}^{2k}\frac{x_i}{a_i}
<1~\mod 2 \bigr\}\\
& -& \#\bigl\{\bfx\in \bbz^{2k+1}
~|~0<x_i<a_i~\hbox{and}~1<\sum_{j=0}^{2k}\frac{x_i}{a_i} <2~\mod 2 \bigr\}.
\end{eqnarray}
Using a formula of Eisenstein, Zagier (cf. \cite{Hir71}) has
rewritten this  as:
\begin{equation}\label{zagfor}
\grt(M^{4k}(\bfa))=\frac{(-1)^k}{N}\sum_{j=0}^{N-1}\cot\frac{\pi(2j+1)}{2N}
\cot\frac{\pi(2j+1)}{2a_1}\cdots \cot\frac{\pi(2j+1)}{2a_{2k+1}},
\end{equation}
where $N$ is any common multiple of the $a_i$'s.

Both formulas are quite well suited to
computer use. A simple C code called {\tt sig.c},
which for any $m$-tuple with $m=2k+1=5,7,9$, computes
the signature $\grt(\bfa):=\grt(M^{4k}_\bfa)$ and the diffeomorphism type of the link  using
either of the above formulas can be found in \cite{BGKT03}.

\begin{exmp}
Let us consider the Brieskorn-Pham link $L(5,3,2,2,2)$.
By Brieskorn Graph Theorem this is a homotopy 7-sphere. One can easily compute
the signature using (2.14) to find out that $\tau(L(5,3,2,2,2))=8$. Hence
 $L(5,3,2,2,2)=\Sigma^7_1$ is an exotic 7-sphere and it is
called Milnor generator (all others can be obtained from it by taking
connected sums). It is interesting to note that one does not need a computer to
find the signature of $L(6k-1,3,2,2,2).$ This was done in Brieskorn original paper \cite{Bri66} 
where he used the combinatorial formula \ref{briescomb} to show that all the 28 diffeomorphism 
types of homotopy 7-spheres are realized by taking $k=1,\cdots,28.$
\end{exmp}

\begin{ques} Suppose
that $\calz_f$ is Fano. How can one prove the existence of
a K\"ahler-Einstein metric on  $\calz_f$? When this can be done successfully we
automatically get a Sasakian-Einstein metric on the link $L_f$.
\end{ques}

\bigskip
\section{Calabi Conjecture I}
\bigskip
Recall that on a K\"ahler manifold the Ricci curvature 2-form $\rho_\omega$
of any K\"ahler metric represents the cohomology class $2\pi c_1(M)$.
The well-known Calabi Conjecture is the question whether or not the
converse is also true. To be more specific we begin with a couple of
definitions

\begin{defn}
Let $(M,J,g,\omega_g)$ be a compact
K\"ahler manifold. The K\"ahler cone of $M$
$$K(M)=\Bigl\{[\omega]\in H^{1,1}(M,\bbc)\cap H^2(M,\bbr)\ \ | \
[\omega]=[\omega_h]\ \qtq{for some K\"ahler metric $h$}\!\!\!\!\!\Bigr\}
$$
is the set of all possible K\"ahler classes on $M$.
\end{defn}

The global $i\partial\bar\partial$-lemma provides for a very simple description of
the space of K\"ahler metrics $\calk_{[\omega]}$. Suppose we have a
K\"ahler metric $g$ in a with K\"ahler class $[\omega_g]=[\omega]\in K(M)$.
If $h\in \calk_{[\omega]}$ is another K\"ahler metric
then, up to a constant, there exists a global function $\phi\in C^\infty(M,\bbr)$
such that $\omega_h -\omega_g=i\partial\bar\partial\phi$. We could fix the
constant by requiring, for example, that $\int_M\phi\omega_g^n=0$. Hence, we have

\begin{cor}
 Let $(M,J,g,\omega_g)$ be a compact K\"ahler manifold
with $[\omega_g]=[\omega]\in K(M)$. Then, relative to the metric $g$
the space $\calk_{[\omega]}$ of all K\"ahler metric in the same K\"ahler class can be described 
as
$$\calk_{[\omega]}=\{\phi\in C^\infty(M,\bbr)\ | \ \omega_h=\omega_g+
i\partial\bar\partial\phi>0, \ \  \int_M\phi\omega_g^n=0\},$$
where the 2-form
$\omega_h>0$ means that $\omega_h(X,JY)$ is a Hermitian metric on $M$.
\end{cor}

We have the following theorem

\begin{thm} Let $(M,J,g,\omega_g)$ be a compact K\"ahler manifold,
$[\omega_g]=\omega\in K(M)$ the corresponding K\"ahler class and $\rho_g$
the Ricci form. Consider any real $(1,1)$-form $\Omega$ on $M$ such that
$[\Omega]=2\pi c_1(M)$. Then there exists a unique K\"ahler metric
$h\in \calk_{[\omega]}$ such that $\Omega=\rho_h$.
\end{thm}

The above statement is the celebrated Calabi Conjecture which
was posed by Eugene Calabi in 1954. The conjecture in its full generality was
eventually proved by Yau in 1978 \cite{Yau78}.

Let us reformulate the problem using the global $i\partial\bar\partial$-lemma.
We start with a given K\"ahler metric $g$ on $M$ in the K\"ahler class $[\omega_g]
=[\omega]$. Since $\rho_g$ also represents $2\pi c_1(M)$ there exists
a globally defined function $f\in C^\infty(M,\bbr)$ such that
$$\rho_g-\Omega=i\partial\bar\partial f.$$
Appropriately, $f$ may be called {\it a discrepancy potential function} for the
Calabi problem
and we could fix the constant by asking that $\int_M(e^f-1)\omega_g=0$.

Now, suppose that the desired solution of the problem is a metric
$h\in\calk_{[\omega]}$. We know that the K\"ahler form of $h$ can be
written as
$$\omega_h=\omega_g+i\partial\bar\partial \phi,$$
for some smooth function $\phi\in C^\infty(M,\bbr)$. We normalize
$\phi$ as in the previous corollary. Combining these two equations we see that
$$\rho_h-\rho_g=i\partial\bar\partial f.$$
If we define a smooth function $F\in C^\infty(M,\bbr)$ relating the
volume forms of the two metrics $\omega_h^n=e^F\omega_g^n$ then the left-hand
side of the above equation takes the following form
$$i\partial\bar\partial F=\rho_h-\rho_g=i\partial\bar\partial f,$$
or simply $i\partial\bar\partial(F-f)=0$. Hence, $F=f+c$. But
since we normalized $\int_M(e^f-1)\omega_g=0$ we must have $c=0$.
Thus, $F=f$, or  $\omega_h^n=e^f\omega_g^n$. We can now give two
more equivalent formulations of the Calabi Problem.

\begin{thm}
Let $(M,J,g,\omega_g)$ be a compact K\"ahler manifold,
$[\omega_g]=[\omega]\in K(M)$ the corresponding K\"ahler class and $\rho_g$
the Ricci form. Consider any real $(1,1)$-form on $M$ such that
$[\Omega]=2\pi c_1(M)$. Let $\rho_g-\Omega=i\partial\bar\partial f,$
with $\int_M(e^f-1)\omega_g=0$.
\begin{enumerate}
\item There exists a unique K\"ahler metric $h\in\calk_{[\omega]}$
whose volume form $\omega_h^n$ equals $e^f\omega_g^n$.
\item Let $(\calu,z_1,\ldots,z_n)$ be a local complex chart on $M$ with
respect to which the metric $g=(g_{i\bar{j}})$. Then, up to a constant,
there exists a unique smooth function $\phi$ in $\calk_{[\omega]}$, which
satisfies the following equation
$$\frac{{\rm det}\Bigl(g_{i\bar{j}}+
\frac{\partial^2\phi}{\partial z_i\partial\bar z_j}\Bigr)}{
{\rm det}(g_{i\bar{j}})}=e^f.$$
\end{enumerate}
\end{thm}

The equation in (2) is called the Monge-Amp\`ere equation.
Part (1) gives a very simple geometric
characterization of the Calabi-Yau theorem. On a compact K\"ahler manifold
one can always find a metric with arbitrarily prescribed volume form.
The uniqueness part of this theorem was already proved by Calabi.
This part involves only Maximum Principle. The existence proof
uses the continuity method discussed briefly in Section 9 and it involves several difficult
a priori estimates. These were found by Yau in 1978. We have the following:

\begin{cor}
Let $(M,J,g,\omega_g)$ be a compact K\"ahler manifold.
\begin{enumerate}
\item If $c_1(M)>0$ then $M$ admits a K\"ahler metric of positive Ricci curvature.
\item If $c_(M)=0$ then $M$ admits a unique K\"ahler Ricci-flat metric.
\end{enumerate}
\end{cor}

It is ``folklore" that the Calabi-Yau Conjecture is also true for
compact orbifolds. In the context of Sasakian geometry with its characteristic                   
foliation, the transverse space $\calz$ is
typically a compact K\"ahler orbifold.
In the context of foliations a {\bf transverse Yau theorem}
was proved by El Kacimi-Alaoui in 1990 \cite{ElK}.

\begin{thm} If $c_1(\calf_\xi)$  is represented
by a real basic $(1,1)$-form $\grr^T,$ then it is the Ricci curvature form
of a unique transverse K\"ahler form $\gro^T$ in the same basic cohomology
class as $d\eta.$
\end{thm}

In the language of
positive Sasakian manifolds this theorem provides the basis for proving \cite{BGN03a}

\begin{thm}\label{posRic}
Any positive Sasakian manifold $(M,g)$ admits
a Sasakian metric $g^\prime$ of positive Ricci curvature.
\end{thm}

There is a similar statement for negative and null Sasakian structures. This is studied in a 
forthcoming article \cite{BGM04}.

\begin{cor} Let $L(f)$ be the link of an isolated hypersurface singularity
where $f$ is weighted homogeneous polynomial of weight $\bfw$ and degree $d$.
If $\sum_iw_i-d>0$ then  $L(f)$ admits a Sasakian metric of positive Ricci curvature.
In particular,
\begin{enumerate}
\item $\#_kM_\infty=\#_k(S^2\times S^3)\simeq L(X_{k+1}\subset \bbp(1,1,1,k))$ admits
Sasakian metric of positive Ricci curvature for all $k\geq1.$
\item The links
$L(6k-1,3,2,\ldots,2)\simeq\Sigma^{4n+3}_k$ and $L(2p+1,2,\ldots,2)\simeq\Sigma^{4n+1}_p$
admit Sasakian metric of positive Ricci curvature for each $n,k$ and $p$. Hence,
all homotopy spheres that bound parallelizable manifolds admit metrics of
positive Ricci curvature.
\end{enumerate}
\end{cor}
Part (1) is the 5-dimensional case of a well known theorem of Sha and Yang \cite{ShYa91}.
Their theorem asserts the existence of positive Ricci curvature metrics on connected sums
of product of spheres. Part (2) is in its final form a theorem of Wraith \cite{Wra97}.
Both papers rely on techniques of surgery theory. The proofs given in \cite{BGN03a, BGN03b} 
are 
completely different being a consequence of the orbifold version of Yau's theorem.

\bigskip
\section{Calabi Conjecture II -- K\"ahler-Einstein Metrics}
\bigskip

As pointed out in the previous section
one of the consequences of the Yau's theorem is that a compact
K\"ahler manifold with $c_1(M)=0$ must admit a Ricci-flat, hence, Einstein
metric. More generally, we can consider existence of K\"ahler-Einstein metrics
with arbitrary Einstein constant $\lambda$. By scaling we can
assume that $\lambda=0,\pm1$. Specifically, let $(M,J,g,\omega_g)$
be a compact K\"aher manifold. We would like to know if one can always find
a K\"ahler-Einstein metric $h\in\calk_{[\omega_g]}$.
Recall that on a K\"ahler-Einstein manifold
$\rho_g=\lambda\omega_g$.  This implies that $2\pi c_1(M)=\lambda [\omega_g]$.
Now, if $c_1(M)>0$ we must have $\lambda=+1$ because $[\omega_g]$ is
the K\"ahler class. Similarly, when $c_1(M)<0$ the only allowable
sign of a K\"ahler-Einstein metric on $M$ is $\lambda=-1$. Clearly,
when $c_1(M)=0$ we must have $\lambda=0$ as $[\omega_g]\not=0$.
As we have already pointed out the $\lambda=0$ case follows from
Yau's solution to the Calabi conjecture. For the reminder of this lecture
we shall assume that $\lambda=\pm1$.

Let $(M,J,g,\omega_g)$ be a K\"ahler manifold and $[\omega_g]=[\omega]\in K(M)$
the K\"ahler class.
Let us reformulate the existence
problem using the global $i\partial\bar\partial$-lemma.
Suppose there exists an Einstein metric
$h\in\calk_{[\omega]}$.
Starting with the original K\"ahler metric $g$ on $M$ we have
a globally defined function $f\in C^\infty(M,\bbr)$ such that
\begin{equation}\label{discrepancy}
\rho_g-\lambda\omega_g=i\partial\bar\partial f.
\end{equation}
As before we will call
$f$ {\it a discrepancy potential function}.
We also fix the constant by asking that $\int_M(e^f-1)\omega_g=0$.
Let $h\in\calk_{[\omega]}$ be an Einstein metric for which
$\rho_h=\lambda\omega_h$. Using the global $i\partial\bar\partial$-lemma
once again we have a globally defined function $\phi\in C^\infty(M,\bbr)$
such that
$\omega_h-\omega_g=i\partial\bar\partial\phi$. We shall fix the
constant in $\phi$ later. Using these two equations we
easily get
$$\rho_g-\rho_h=i\partial\bar\partial(f-\lambda\phi).$$
Defining $F$ so that $\omega_h^n=e^F\omega_g^n$ we can write this
equation as
$$i\partial\bar\partial F= i\partial\bar\partial(f-\lambda\phi).$$
This implies that $F=f-\lambda\phi+c$. We have already fixed
the constant in $f$ so $c$ depends only on the choice of $\phi$.
We can make $c=0$ by choosing $\phi$ such that
$\int_M(e^{f-\lambda\phi}-1)\omega_g^n=0$.
Hence, we have the following

\begin{prop} Let $(M,J)$ be a compact K\"ahler manifold
with $\lambda c_1(M)>0$, where $\lambda=\pm1$. Let $[\omega]\in K(M)$
be a K\"ahler class and  $g,h$ two K\"ahler metrics in $\calk_{[\omega]}$
with Ricci forms $\rho_g,\rho_h$. Let $f,\phi\in C^\infty(M,\bbr)$
be defined by $\rho_g-\lambda\omega_g=i\partial\bar\partial f,$
$\omega_h-\omega_g=i\partial\bar\partial \phi$. Fix the relative
constant of $f-\lambda\phi$ by setting
$\int_M(e^{f-\lambda\phi}-1)\omega_g=0$. Then the metric $h$ is Einstein
with Einstein constant $\lambda$ if and only if $\phi$ satisfies
the following Monge-Amp\`ere equation
$$\omega_h^n=e^{f-\lambda\phi}\omega_g^n,$$
which in a local complex chart $(\calu,z_1,\ldots,z_n)$ written as
$$\frac{{\rm det}\Bigl(g_{i\bar{j}}+
\frac{\partial^2\phi}{\partial z_i\partial\bar z_j}\Bigr)
}{{\rm det}(g_{i\bar{j}})}=e^{f-\lambda \phi}.$$
\end{prop}

Note that by setting $\lambda=0$ we get the Monge-Amp\`ere equation
for the original Calabi problem.
The character of the Monge-Amp\`ere equation  above very much depends
on the choice of $\lambda$. The case of $\lambda=-1$ is actually
the simplest as the necessary a priori $C^0$-estimates can be
derived using the Maximum Principle. This was done by Aubin \cite{Aub76} and independently
by Yau \cite{Yau78}. We have

\begin{thm} Let $(M,J,g,\omega_g)$ be a compact K\"ahler
manifold with $c_1(M)<0$. Then there exists a unique K\"ahler metric
$h\in\calk_{[\omega_g]}$ such that $\rho_h=-\omega_h$.
\end{thm}

When $k=+1$ the problem is much harder. It has been known for quite some time
that there are non-trivial obstructions to the existence of K\"ahler-Einstein metrics.
Let $\gh(M)$ be the complex Lie algebra of all holomorphic vector fields
on $M$. Matsushima \cite{Mat57} proved that on a compact K\"ahler-Einstein
manifold with $c_1(M)>0,$ ${\gh}(M)$ must be reductive, i.e.
${\gh}(M)=Z({\gh}(M))\oplus[{\gh}(M),{\gh}(M)]$ where $Z({\gh}(M))$ denotes the center of 
$\gh(M).$ 
Now, suppose $(M,g,J,\omega)$ is a K\"ahler manifold and let
$f$ be the discrepancy potential defined by (\ref{discrepancy}). Further, let $X\in {\gh}(M)$ and
define
\begin{equation}
F(X)=\int_MX(f)d{\rm vol}_g.
\end{equation}
At first glance $F(X)$ appears to depend on the K\"ahler metric. However,
Futaki shows that this is not the case:
$F(X)$ does not depend on the choice of the metric in $h\in\calk_{[\omega_g]}$.
Hence, $F:\gh(M)\ra{1.2}\bbc$ is well-defined and it is called
the {\it Futaki functional} or {\it character} \cite{Fut83,Fut87, Fut88}.
In particular, if an Einstein metric
exists than we can choose $f$ constant. Hence,

\begin{cor}
Suppose $(M,g,J,\omega)$ admits a K\"ahler-Einstein metric $h\in\calk_{[\omega_g]}$.
Then $F$ must be identically zero.
\end{cor}

\begin{rem}
Futaki also showed that there are Fano manifolds for which ${\gh}(M)$ is reductive, but
$F$ is non-trivial \cite{Fut88}. A folklore conjecture attributed to Calabi asserted that
in the case $\gh(M)=0$ there are no obstructions to finding a K\"ahler-Einstein metric.
This conjecture was disproved by Tian. First, Ding and Tian \cite{DiTi92} constructed
an example of an orbifold del Pezzo surface with $\gh(M)=0$ and no K\"ahler-Einstein metric.
Later Tian found an example of a smooth Fano 3-fold with $\gh(M)=0$ and no
K\"ahler-Einstein metric
\cite{Tia99}. In \cite{Tia97} Tian shows that two different conditions are
necessary for the existence of a K\"ahler-Einstein metric. One condition involves the
generalized Futaki functional of every special degeneration of the manifold. The other
condition is Mumford stability with respect to a certain polarization. Tian
conjectures that the two conditions are equivalent and that they are also
sufficient. This conjecture is still open. However, even if it is true neither of the
conditions are easily checked for an arbitrary compact Fano manifold (orbifold).
In principle, one should be able to compute generalized Futaki invariants for any Fano
hypersurface
$X_d\subset\bbp(\bfw)$.
(see, for example, Lu \cite{Lu99}, and Yotov \cite{Yot99} for
the computation of generalized Futaki invariants in the case of smooth complete intersections).
\end{rem}
\bigskip
\section{The Continuity Method and K\"ahler-Einstein Orbifolds}
\bigskip

Let us briefly describe the main aspects of the continuity method. Let's say we are trying
to show existence of K\"ahler-Einstein metric of positive sign.
Here one tries to solve the Monge-Amp\`ere equation
$$\frac{\det(g_{i\barj}+\partial_i\bar{\partial}_{\barj}\phi_t)}{
\det(g_{i\barj})}= e^{-t\phi_t+f}, \quad
g_{i\barj}+\partial_i\bar{\partial}_{\barj}\phi_t>0$$
for $t\in [0,1].$  Yau's Theorem tells us that this has a solution for $t=0,$
and we try to solve this for $t=1$, where the metric will be K\"ahler-Einstein.
The so called continuity method sets out to show that the interval where solutions exist
is both open and closed. Openness follows from the Implicit Function Theorem,
but there are well known obstructions to closedness. This problem has been
studied most recently by Demailly
and Koll\'ar who work in the orbifold category \cite{DeKo01}. Closedness is equivalent to
the uniform boundedness of the integrals
$$\int_{\calz}e^{-\grg t\phi_t}\gro_0^n$$
for any $\grg\in (\frac{n}{n+1},1),$ where $\gro_0$ is the
K\"ahler form of $h_0.$ This means that the {\it multiplier ideal sheaf} of Nadel \cite{Nad90}
$\calj(\grg\phi)=\calo_\calz$ for all $\grg\in (\frac{n}{n+1},1).$

We will illustrate how the method works for links.
The approach developed by Demailly and Koll\'ar in \cite{DeKo01}
yields the following general theorem:

\begin{thm}\label{DKgeneral}
Let $X^{orb}$ be a compact, $n$-dimensional orbifold such that
 $K_{X^{orb}}^{-1}$ is ample.
The continuity method produces a
\Ke metric on $X^{orb}$ if the following holds:
There is a $\gamma>\tfrac{n}{n+1}$ such that
 for every
$s\geq 1$ and for every holomorphic section
$\tau_s\in H^0(X^{orb}, K_{X^{orb}}^{-s})$
the following integral is finite:
$$
\int  |\tau_s|^{-\frac{2\gamma}{s}}\omega_0^n\  <\  +\infty.
$$
\end{thm}

In general, this condition is not hard to check. For hypersurfaces
the situation is somewhat simpler and one gets

\begin{cor}\label{DKhypersurface} Let $\calz_f=:X_d\subset\bbp(w)$ be a hypersurface
in $\bbp(w)$ given by the vanishing of the weighted homogeneous polynomial $f$ of
weight $\bfw$ and degree $d$. Let $Y_f:=\{\{f=0\}\subset \bbc^{n+1}\}$ so that
$\calz_f=(Y_f\setminus\{0\})/\bbc_\bfw^*$. Assume $\calz_f$ is Fano, that is
$d<\sum w_i$.
The continuity method produces a
K\"ahler-Einstein metric on $\calz_f$ if the following holds:
There is a $\gamma>\tfrac{n}{n+1}$ such that
for every weighted homogeneous polynomial $g$ of
weighted degree $s(\sum w_i-d))$, not identically zero on $Y_f$, the function
$$
|g|^{-\gamma/s}\qtq{is locally $L^2$ on $Y_f\setminus\{0\}$.}
$$
\end{cor}

\begin{exmp}\label{LDP} [{\bf Log Del Pezzo Surfaces}]
One can take $X_d\subset \bbp(w_0,w_1,w_2,w_3)$ with $I=\hbox{Ind}(X_d)=\sum_iw_i-d>0$.
Assuming $X_d$ has only isolated orbifold singularities one can classify all such
log del Pezzo surfaces and check if the conditions of Corollary \ref{DKhypersurface}
are satisfied. This was done by Johnson and Koll\'ar in \cite{jk1} when $I=1$. In some
cases the existence question was left open and more recently Araujo completely finished
the analysis in \cite{Ara02}. On the other hand, the similar question can
be considered for an arbitrary index $I\geq1$. This was done to a limited extend in
\cite{BGN03c}.
That is, we enumerated all log del Pezzo surfaces for $2<I\leq10$ which can possibly admit
K\"ahler-Einstein
metrics as a consequence of the Corollary \ref{DKhypersurface}. Sometimes we were able
to prove the existence, but unlike in the $I=1$ situation, it has not been done for all
the candidates. Furthermore, it remains to show that for $I>10$ there are no
examples of log del Pezzo surfaces satisfying condition of Corollary \ref{DKhypersurface}.
As a result of this analysis we got Sasakian-Einstein structures on certain
connected sums of $S^2\times S^3$. The table below summarizes the results of
\cite{BGN03c, BGN02b, BG03}.

\vskip 12pt
\begin{center}
\begin{tabular}{|c|c|c|}\hline
$\#_k(S^2\times S^3)$ & {\bf N}=$(n_0,n_1,\ldots,n_7,n_8)$&Example\\
\hline\hline
$k=1$ &$(14+1,0,0,0,0,0,0,0,0)$&$X_{76}\subset\bbp(11,13,21,38)$   \\ \hline
$k=2$ &$(21,2,0,0,0,0,0,0,0)$&$X_{57}\subset\bbp(7,8,19,25)$   \\ \hline
$k=3$ &$(\aleph_0+1,4,2,0,0,0,0,0,0)$&$X_{64}\subset\bbp(7,8,19,32)$  \\ \hline
$k=4$ &$(2\aleph_0+1,0,1,2,0,0,0,0,0)$&$X_{20}\subset\bbp(3,4,5,10)$   \\ \hline
$k=5$ &$(2\aleph_0,0+1,0,1,1,0,0,0,0)$&$X_{28}\subset\bbp(3,5,7,14)$   \\ \hline
$k=6$ &$(2\aleph_0,0,0+1,3,0,5,0,0,0)$&$X_{18}\subset\bbp(2,3,5,9)$  \\ \hline
$k=7$ &$(0,0,0,0+1,0,\aleph_0,0,0,0)$&$X_{8k+4}\subset\bbp(2,2k+1,2k+1,4k+1)$  \\ \hline
$k=8$ &$(0,0,0,0,0+1,0,0,0,2)$&$X_{10}\subset\bbp(1,2,3,5)$ \\ \hline
$k=9$ &$(0,0,0,0,0,0,0,0,1)$&$X_{16}\subset\bbp(1,3,5,8)$ \\ \hline
\end{tabular}
\end{center}
\vskip 12pt
For each $\#_k(S^2\times S^3)$, with $1\leq k\leq9$ we list
${\bf N}=(n_0,n_1,n_2,\ldots,n_8)$, where $n_i$ is the number of distinct families of links with
complex dimension of deformation parameters
equal to $i$. The largest family constructed this way had complex dimension 8
so that $n_i=0$ when $i>8$ for all $k=1,\ldots,9$. Also, the method produced no
examples for $k>9$. Furthermore, for instance, $n_0=2\aleph_0$ means that
there are two distinct infinite sequence of examples which have no deformation
parameters. We include regular examples in the count by writing ``$+1$" where
appropriate. In the third column we give
an example with the largest moduli.
\end{exmp}

\begin{rem} Example \ref{LDP} forces an obvious question. Are there Sasakian-Einstein
structures on $\#_k(S^2\times S^3)$ for arbitrary $k$? Recently Koll\'ar \cite{Kol04}
has been able to answer this question in the affirmative. His method differs substantially
from the one described here. The idea is to consider Seifert bundles over smooth
surfaces, but with a non-trivial orbifold structure. Such a construction is more flexible in obtaining 
log del Pezzo surfaces with orbifold K\"ahler-Einstein metrics and
more complicated topology. In particular, Koll\'ar proves

\begin{thm} For every integer 
$k\geq6$ there are infinitely many complex $(k-1)$-dimensional
families of Einstein metrics on $\#_k(S^2\times S^3)$.
\end{thm}

Combining this remarkable result with the links of hypersurfaces in Example \ref{LDP},
we get the following

\begin{cor} Let $M$ be any compact, smooth, simply-connected 5-manifold which is spin and
has no torsion in $H_2(M,\bbz)$. Then $M$ admits a Sasakian-Einstein metric.
\end{cor}
\end{rem}

\begin{exmp}\label{BPlink} [{\bf Brieskorn-Pham Links}] Now, we
consider Brieskorn--Pham links as defined in Equation (\ref{BP}). Let
$Y(\bfa):=(\sum_{i=0}^nz_i^{a_i}=0)\subset \bbc^{n+1}.$
One can easily see that $d={\rm lcm}(a_i:i=0,\dots,n)$ is the degree of $f_\bfa(\bfz)$ and
$w_i=d/a_i$ are the weights.
The transverse space
$\calz(\bfa):=(Y(\bfa)\setminus\{0\})/\c^*$  is a Fano orbifold if and only if
$$1<\sum_{i=0}^n\frac{1}{a_i}.$$
\noindent
More generally, we consider weighted homogeneous perturbations
$$
Y(\bfa,p):=(\sum_{i=0}^nz_i^{a_i}+p(z_0,\dots,z_n)=0)\subset \bbc^{n+1},$$
where weighted ${\rm degree}(p)=d$.
The genericity
condition we need, which is always satisfied by $p\equiv 0$ is:
The intersections of $Y(\bfa,p)$ with any number of hyperplanes
$(z_i=0)$
are all smooth outside the origin.

The continuity methods produces the following sufficient conditions
for the quotient 
$Y(\bfa,p)/\bbc^*$ to admit a K\"ahler-Einstein metric \cite{BGK03}: 

\begin{thm}\label{spheres} Let
$\calz(\bfa,p)$ be the transverse space of
a perturbed Brieskorn-Pham link $L(\bfa,p)$.
Let $C_i={\rm lcm}(a_0,\ldots,\hat a_i,\ldots,a_n)$, $b_i={\rm gcd}(C_i,a_i)$.
Then $\calz(\bfa,p)=Y(\bfa,p)/\bbc^*$ is Fano and it
has a K\"ahler-Einstein metric if
\begin{enumerate}
\item
$1<\sum_{i=0}^n \frac{1}{a_i}$,
\item
$\sum_{i=0}^n \frac{1}{a_i}<1+\frac{n}{n-1}\min_{i}\{\frac{1}{a_i}\}$, and
\item
$\sum_{i=0}^n \frac{1}{a_i}<1+\frac{n}{n-1}\min_{i,j}\{\frac{1}{b_ib_j}\}$.
\end{enumerate}
In this case the link $L(\bfa,p)$ admits a Sasakian-Einstein metric
with one-dimensional isometry group.
\end{thm}
\end{exmp}

\bigskip
\section{Sasakian-Einstein Structures on Brieskorn-Pham Links}
\bigskip

In this section we will discuss some consequences of Theorem \ref{spheres} of the
previous section. We will investigate two separate cases: when $L(\bfa)$ is a homotopy sphere
and (2) when $L(\bfa)$ is a rational homology sphere with non-vanishing torsion.
If $L(\bfa)$ is a homotopy sphere, for a fixed $n$,
there are only finitely many examples of $\bfa$'s satisfying all three conditions of
Theorem \ref{spheres}. However, the number of examples as well as the moduli
grows doubly exponentially with each odd dimension. One can list all
solution in dimensions 5 and 7 without difficulties. However, already in
dimension 9 that task is too overwhelming. It is quite clear that
many of our links will actually be diffeomorphic to standard spheres.
Hence, let us begin with a remark concerning what is known about Einstein metrics on
spheres in general.
\begin{rem}
Any standard sphere $S^n$, $n>1$, admits a metric of constant positive sectional
curvature. These canonical metrics are $SO(n+1)$-homogeneous and Einstein, i.e.,
the Ricci curvature tensor is a constant positive multiple of the metric.
The spheres $S^{4m+3}$, $m>1$ are known to have
another $Sp(m+1)$-homogeneous Einstein metric discovered by
Jensen \cite{Jen73}. The metric is obtained from the ``quaternionic Hopf fibration"
$S^3\rightarrow S^{4m+3}\rightarrow \bbh\bbp^m$. Since both base
and fiber are Einstein spaces with positive Einstein constant we
obtain two Einstein metrics in the canonical variation. The second metric
is also called ``squashed sphere" metric in some physics literature.
In addition, $S^{15}$ has a third ${\rm  Spin}(9)$-invariant
homogeneous Einstein metric discovered by Bourguignon and Karcher in 1978 \cite{BuKa78}.
The existence of such a metric has to do with the fact that $S^{15}$, in addition
to fibering over $\bbh\bbp^2$, also fibers over $S^8$ with fiber $S^7$.
Thus the15-sphere admits 3 different homogeneous Einstein metrics.
Ziller proved that these are the only
homogeneous Einstein metrics on spheres \cite{Zil82}.
B\"ohm obtained infinite sequences of
non-isometric Einstein metrics, of positive scalar curvature, on
$S^5$, $S^6$, $S^7$, $S^8$, and $S^9$ \cite{Boe98}. B\"ohm's metrics are
of cohomogeneity one and they are not only the first
inhomogeneous Einstein metrics on spheres but also the first
non-canonical Einstein metrics on even-dimensional spheres.
\end{rem}

\begin{exmp} \label{5s} [{\bf Sasakian-Einstein Metrics on} $S^5$]
Consider $L(2,3,7,m)$. These are homotopy spheres as long as $m$ is relatively prime to
at least two of $2,3,7$. It is easy to see that $L(2,3,7,m)$ satisfies the condition
of Theorem \ref{spheres}  if $5\leq m\leq 41$ which gives $27$ cases.
The link $L(2,3,7,35)$ admits deformations, i.e.,
$C(u,v)$ is any sufficiently general homogeneous septic polynomial, then
the link of
$$
z_0^2+x_1^3+C(z_2,z_3^5)
$$
also gives a Sasakian-Einstein metric
on $S^5$.  The relevant automorphism group of $\bbc^4$ is
$$
(z_0,z_1,z_2,z_3)\mapsto
(z_0,z_1,\alpha_2z_2+\beta z_3^5,\alpha_3z_3).
$$
Hence we get a $2(8-3)=10$ real dimensional family of Sasakian-Einstein metrics on $S^5$.
There are other examples, 68 in total, and we get \cite{BGK03}

\begin{thm}
On $S^5$ there are at least 68 inequivalent families of
Sasakian-Einstein metrics. Some
of these families admit  non-trivial continuous Sasakian-Einstein deformations.
The biggest constructed family has has real dimension 10.
\end{thm}
\end{exmp}
\bigskip
\begin{exmp} \label{7s}  [{\bf Sasakian-Einstein Metrics on Homotopy 7-Spheres}]
Similarly $L(2,3,7,43,43\cdot31)$ is the standard 7-sphere
with a $2(43-2)=82$-dimensional family of Sasakian-Einstein metrics.
One can do a computer search of all homotopy 7-spheres that satisfy the numerical
conditions of Theorem \ref{spheres}. One finds 8610 such links. An additional
computation of the Hirzebruch signature of the parallelizable manifold whose boundary is $L$ 
shows that they are more or less evenly
distributed among the 28 oriented diffeomorphism classes. This way we get 
\cite{BGK03,BGKT03}

\begin{thm}
Let $\Sigma^7_i$, be a homotopy 7-sphere corresponding
to the element $i\in bP_8\simeq\bbz_{28}\simeq \Theta_7$ in the Kervaire-Milnor group.
$\Sigma_i^7$ admits at least
$n_i$ inequivalent families of
Sasakian-Einstein  metrics,
where $(n_1,\ldots,n_{28})=${\rm (376, 336, 260, 294, 231,
284, 322, 402, 317, 309, 252, 304, 258,
390,409, 352, 226, 260, 243, 309, 292, 452, 307, 298, 230, 307, 264, 353)}, giving a total
of $8610$ cases. In each oriented diffeomorphism class some of the families depend on a
moduli.
In particular, the standard 7-sphere $\Sigma_{28}^7$  admits an 82-dimensional family of
inequivalent Sasakian-Einstein metrics.
\end{thm}

\end{exmp}

\begin{exmp}  [{\bf Sasakian-Einstein Metrics on Kervaire Spheres}]
\noindent
Let $\{c_i\}$ be an infinite sequence defined by the recursion relation
$$c_{k+1}=c_1\cdots c_k+1=
c_k^2-c_k+1,\ \ \ \ \ c_1=2.$$
Consider sequences of the form
$L(\bfa)=L(2c_1,2c_2,\dots,2c_{m-2},2, a_{m-1})$
where $a_{m-1}$ is relatively prime to all the other $a_i$s. Such $L(\bfa)$ are rational
homotopy spheres.
The condition of Theorem~\ref{spheres}
is satisfied if $2c_{m-2}<a_{m-1}<2c_{m-1}-2$.
In particular, we can ask for  $a_{m-1}$ to be prime and estimate the number of
primes in the range $(2c_{m-2}, 2c_{m-1}-2)$ which gives double exponential growth in $m$
by the Prime Number Theorem. For odd $m$,  $L(\bfa)$  the standard sphere if one of the 
exponents of $\bfa$
equals to $\pm1 \mod 8$ and it is the Kervaire sphere if one of the exponents equals $\pm 3
\mod 8$
\cite{Bri66}.
It is easy to check for all values of $m$ that we get at least one solution
of both types. Hence, we get \cite{BGK03}
\begin{prop}\label{Ke} Theorem~\ref{spheres} yields a doubly exponential number of
inequivalent Sasakian-Einstein metrics on both the standard and Kervaire spheres in every
odd dimension $4m+1$.
\end{prop}
\end{exmp} 

\begin{exmp}
$L(\bfa)=L(2,3,7,43,1807, 3263443, 10650056950807, m)$
is just the standard 13-sphere for any suitably chosen $m$ as $bP_{14}$ is trivial.
If we choose $m=(10650056950807-2)\cdot 10650056950807$ we get
$2(10650056950807-2)$-dimensional family of deformations.
By contrast the only Einstein metric on $S^{13}$ known previously was the
the canonical one.
\end{exmp}

All these examples point towards the following:
\begin{conj} All odd-dimensional spheres that bound parallazllizable manifolds
admit Sasakian-Einstein metrics.
\end{conj}

The conjecture is true in dimension $4m+1$ by Proposition~\ref{Ke}. It is also true
in dimension 7. In addition, using computer programs we were able to verify
that the conjecture holds in dimensions 11 and 15. Computational verification
in arbitrary dimension $4m+3$ is not possible. On the other hand
it does appear that Brieskorn-Pham links satisfying the conditions of
Theorem~\ref{spheres} realize all oriented diffeomorphism types of
homotopy spheres in every dimension.

\begin{exmp} \label{rhs} [{\bf S-E Structures on Rational Homology Spheres}]
Our final examples of Brieskorn-Pham links is that of
$L(m,m,\ldots,m,k)$ with ${\rm gcd}(k,m)=1$. By Brieskorn Graph Theorem this
is a rational homology sphere in every dimension. The
conditions of Theorem~\ref{spheres} are satisfied as long as $k>m(m-1)$.
The homology of $L(m,m,\ldots,m,k)$ contains torsion in $H_{m-1}(L,\bbz)$. It's order can be
easily
computed and it is $k^{b_{m-1}}$, where $b_{m-1}$ is the $(m-1)^{th}$ Betti number
of the null link $L(m,\ldots,m)$ which is a regular circle bundle over
Fermat hypersurface of degree $m$ in $\bbp^{m-1}$. For example, with the appropriate
restriction on $k$ we have
$|H_2(3,3,3,k)|=k^2$ and $|H_3(4,4,4,4,k)|=k^{21}$, and so on.  In particular, we get 
\cite{BG03p}
\begin{prop} In each odd dimension greater than 3 there are infinitely many
smooth, compact, simply-connected rational homology spheres admitting Sasakian-Einstein
structures. 
\end{prop}

Torsion groups of each of the links  $L(m,m,\ldots,m,k)$ can also be computed
using an algorithm conjectured by Orlik \cite{Or72} and proved by 
Randell in some special cases \cite{Ran75}. In particular,
Orlik's conjecture is true for all Brieskorn-Pham links and can be used to
compute torsion of various examples discussed here
\cite{BG04p}.
Let us consider the 5-dimensional case in more detail. Using 
Theorem~\ref{spheres} and Orlik's algorithm we
get the following list \cite{OrWa75, Or72}
\begin{center}
\begin{tabular}{|c|c|}\hline
$L(\bfa)$ & Torsion \\
\hline\hline
$L(3,3,3,k)$,  ${\rm gcd}(k,3)=1, k>5$&$\bbz_k\oplus\bbz_k$ \\
\hline
$L(2,4,4,k)$, ${\rm gcd}(k,2)=1$, $k>10$&$\bbz_k\oplus\bbz_k$ \\
\hline
$L(2,3,6,k)$,  ${\rm gcd}(k,6)=1$, $k>12$&$\bbz_k\oplus\bbz_k$ \\
\hline
\end{tabular}
\end{center}
The three series above satisfy
$\sum_{i=0}^2\frac{1}{a_i}=0$. In the case when
$\sum_{i=0}^2\frac{1}{a_i}<0$ one can easily see that there are 16
more rational homology 5-spheres which satisfy inequalities of
Theorem~\ref{spheres}. An example of such a link is $L(3,4,4,4)$
whose  2-torsion equals $\bbz_3\oplus\bbz_3\oplus \bbz_3\oplus\bbz_3\oplus
\bbz_3\oplus\bbz_3$. Hence, $L(3,4,4,4)$ is diffeomorphic to
$M_3\#M_3\#M_3$. For the torsion computation as well as the full list we
refer interested readers to \cite{BG04p}.
In particular, we get the following
\begin{thm} The Barden manifold $M_k$ admits Sasakian-Einstein
structure for each $k>5$ prime to 3 and for each $k>10$ prime to 2.
\end{thm}
\end{exmp}

\begin{exmp}
Finally, note that the links in the last table have companions
with non-trivial second Betti number and by Theorem~\ref{spheres} they too admit
Sasakian-Einstein metrics. We list the relevant information
in the table below:
\bigskip
\begin{center}
\begin{tabular}{|c|c|}\hline
$L(\bfa)$ & $b_2(L(\bfa))$ \\
\hline\hline
$L(3,3,3,3n)$, $n>2$&6 \\
\hline
$L(2,4,4,2n)$, ${\rm gcd}(n,2)=1$, $n>5$&3 \\
\hline
$L(2,4,4,4n)$, $n>2$&7 \\
\hline
$L(2,3,6,2n)$,  ${\rm gcd}(n,3)=1$, $n>12$&2 \\
\hline
$L(2,3,6,3n)$,   ${\rm gcd}(n,2)=1$, $n>12$ &4 \\
\hline
$L(2,3,6,6n)$, $n>4$&8 \\
\hline
\end{tabular}
\end{center}
\bigskip
All the above links have 2-torsion equal to $\bbz_n\oplus\bbz_n$ which can be verified
by Orlik's algorithm. In addition,
just as in the case of rational homology 5-spheres, one can see
that there are 16 exceptional cases of links which satisfy the
inequalities of Theorem~\ref{spheres} and have non-vanishing
second Betti number. An example of such a link is
$L(2,4,6,10)$ which has $b_2=1$. Each line of the previous table gives
infinite series of decomposable Barden manifolds of mixed type (i.e., having both a
free part and a torsion in its second homology). For instance, we can rephrase first line
as
\begin{prop}
The manifolds $6M_\infty\#M_n$ admit families of Sasakian-Einstein
structures for any $n>2$. 
\end{prop}

\end{exmp}
\begin{ques} We conclude by asking some questions about Sasakian-Einstein structures on
certain Barden manifolds:
\begin{enumerate}
\item Does every Barden manifold $M_k$ admit a Sasakian-Einstein structure? Is it possible
that certain torsion in $H_2(M,\bbz)$ obstructs the existence of Sasakian-Einstein, or even 
positive Sasakian structures on $M$?
\item Which Barden manifolds with $b_2(M)\not=0$ and non-vanishing 2-torsion
admit Sasakian-Einstein structures?
\end{enumerate}
\end{ques}

\begin{rem}\label{5B} [{\bf Einstein Metrics on Barden Manifolds}]
The following table summarizes what we know about existence of
Sasakian ($\cals$), negative Sasakian  ($\cals_-$),
null  Sasakian ($\cals_0$), positive  Sasakian ($\cals_+$), regular
Sasakian-Einstein ($\cals^{r}$) and 
non-regular Sasakian-Einstein ($\cals^{nr}$) structures
on some Barden manifolds. The last column (OE) indicates if an Einstein
metric other than Sasakian-Einstein is known. Finally, ``some $k$"
means that we know existence of a given structure for some
(possibly infinitely many) $k$'s but we do not know if it exists for
all $k$.
\vskip 12pt
\begin{center}
\begin{tabular}{|c|c|c|c|c|c|c|c|}\hline
$M$ &$\cals$ &$\cals_-$ &$\cals_0$
&$\cals_+$ &$\cals\cale^{r}$ &$\cals\cale^{nr}$ & OE \\
\hline\hline
$X_{-1}$ &?&no&no&no&no&no&yes \\ \hline
$X_{n}$, $n\geq1$  &?&no&no&no&no&no&? \\ \hline
$X_{\infty}$ &yes&no&no&no&no&no&yes \\ \hline
$X_{0}\simeq S^5$  &yes&yes&no&yes&yes&yes&yes \\ \hline
$M_\infty\simeq S^2\times S^3$ &yes&yes&no&yes&yes&yes&yes \\ \hline
$M_{n}, n\not=3k$  &yes&?&?&yes&no&$n>5$&? \\ \hline
$M_{n}, n\not=2k$  &yes&?&?&yes&no&$n>10$&? \\ \hline
$6M_\infty\#M_{n}$  &yes&?&?&yes&no&$n>2$&? \\ \hline
$kM_\infty, 1<k$&yes&some $k$&some $k$&yes&$2\neq k<9$&yes&? \\ \hline
\end{tabular}
\end{center}
\vskip 12pt
The non Sasakian-Einstein metric on Barden manifolds are the following:
$X_{-1}$ is a symmetric space and the metric is Einstein. $S^2\times S^3$
is well known to have infinitely many inequivalent homogeneous Einstein metrics
discovered by Wang and Ziller \cite{WaZi90}. $S^5$ and $S^2\times S^3$ have
infinitely many inequivalent Einstein metrics of cohomogeneity one discovered by
B\"ohm \cite{Boe98}. Finally $X_\infty$ has infinitely many Einstein metrics
recently constructed by several physicists \cite{HSY04, LPP04}.
\end{rem}
\bigskip\bigskip


\def\cprime{$'$} \def\cprime{$'$}
\providecommand{\bysame}{\leavevmode\hbox to3em{\hrulefill}\thinspace}
\providecommand{\MR}{\relax\ifhmode\unskip\space\fi MR }
\providecommand{\MRhref}[2]{%
  \href{http://www.ams.org/mathscinet-getitem?mr=#1}{#2}
}
\providecommand{\href}[2]{#2}

\end{document}